\documentclass[sn-mathphys]{sn-jnl}

\jyear{2022}%

\usepackage[utf8]{inputenc}
\usepackage{microtype}
\usepackage[T1]{fontenc}
\usepackage{dsfont}
\usepackage{amssymb,amsmath,bm,bbm,upgreek}
\usepackage{algorithm,algorithmicx,algpseudocode}
\usepackage{enumitem}
\usepackage{marginnote}
\usepackage[mathscr]{eucal}
\usepackage{graphicx}
\usepackage{xcolor}
\usepackage{booktabs}
\usepackage{subfigure}
\usepackage{dcolumn}
\newcolumntype{d}[1]{D{.}{.}{#1} }
\newcolumntype{e}[1]{D{;}{\,}{#1} }
\newcolumntype{t}[1]{D{;}{\, \cdot \, }{#1} }
\usepackage{tabularx}
\usepackage{array} 
\usepackage{comment}

\newcommand{\RB}{{\mathsf{rb}}}

\newcommand{\bx}{{\bm x}}

\newcommand{\Vs}{{V}}
\newcommand{\Hs}{{H}}

\newcommand{\Hm}{{\mathscr H}}
\newcommand{\Xm}{{\mathscr X}}
\newcommand{\Ym}{{\mathscr Y}}

\newcommand{\Pmad}{{\mathscr P_\mathsf{ad}}}

\newcommand{\calB}{\mathcal{B}}

\newcommand{\calF}{\mathcal{F}}
\newcommand{\calR}{\mathcal{R}}

\newcommand{\calP}{\mathcal{P}}
\newcommand{\calA}{\mathcal{A}}

\newcommand{\calN}{\mathcal{N}}

\newcommand{\trial}{{\mathscr{X}}}
\newcommand{\test}{{\mathscr{Y}}}
\newcommand{\train}{\mathscr{P}_{\mathrm{train}}}
\newcommand{\trainDEIM}{\mathscr{P}_{\mathrm{train}}^L}
\newcommand{\FEspace}{{{V}_{h}}}
\newcommand{\FEtimeL}{{{S}_{\Delta t}}}
\newcommand{\FEtimeC}{{{Q}_{\Delta t}}}
\newcommand{\RBspace}{{{V}_{\ell}}}

\newcommand{\Totaltol}{ {\varepsilon_{\mathrm{tol}} }}
\newcommand{\RBtol}{ {\varepsilon_{\mathrm{tol}}^\RB} }
\newcommand{\DEIMtol}{ {\varepsilon_{\mathrm{tol}}^L} }

\newcommand{\Bmat}{{ \mathsf{B} }}
\newcommand{\Fmat}{{ \mathsf{F} }}
\newcommand{\Pmat}{{ \mathsf{P} }}
\newcommand{\Phimat}{{ \mathsf{\Phi} }}
\newcommand{\Imat}{{ \mathsf{I} }}

\newcommand{\Mspace}{ {\mathsf{M}_h^{\text{space}} }}
\newcommand{\MspaceL}{ {\mathsf{M}_\ell^{\text{space}} }}
\newcommand{\Vspace}{ {\mathsf{V}_h^{\text{space}} }}
\newcommand{\VspaceL}{ {\mathsf{V}_\ell^{\text{space}} }}
\newcommand{\Mlump}{ {\tilde{\mathsf{M}}_h^{\text{space}} }}
\newcommand{\Mtime}{ {\mathsf{M}_h^{\text{time}} }}
\newcommand{\Ntime}{ {\mathsf{N}_h^{\text{time}} }}
\newcommand{\RBmat}{{ \mathsf{\Psi}_\ell }}
\newcommand{\RBmatT}{{ \mathsf{\Psi}_\ell^\top }}

\newcommand{\yfec}{{\mathsf{y}}}

\newcommand{\POD}{\text{POD}}

\theoremstyle{thmstyleone}
\newtheorem{theorem}{Theorem}
\newtheorem{proposition}[theorem]{Proposition}
\newtheorem{lemma}[theorem]{Lemma}
\newtheorem{corollary}[theorem]{Corollary}

\theoremstyle{thmstyletwo}%
\newtheorem{example}{Example}%
\newtheorem{remark}{Remark}%

\theoremstyle{thmstylethree}%
\newtheorem{assumption}{Assumption}%

\DeclareMathOperator*{\argmax}{arg\,max}

\begin{document}

\title[An adaptive certified space-time RB method for nonsmooth PDEs]{An adaptive certified space-time reduced basis method for nonsmooth parabolic partial differential equations}

\author*[1]{\fnm{Marco} \sur{Bernreuther}}\email{marco.bernreuther@uni-konstanz.de}

\author[1]{\fnm{Stefan} \sur{Volkwein}}\email{stefan.volkwein@uni-konstanz.de}

\affil[1]{\orgdiv{Department of Mathematics and Statistics}, \orgname{University of Konstanz}, \postcode{78457}~\city{Konstanz}, \country{Germany}}

\abstract{In this paper, a nonsmooth semilinear parabolic partial differential equation (PDE) is considered. For a reduced basis (RB) approach, a space-time formulation is used to develop a certified a-posteriori error estimator. This error estimator is adopted to the presence of the discrete empirical interpolation method (DEIM) as approximation technique for the nonsmoothness. The separability of the estimated error into an RB and a DEIM part then guides the development of an adaptive RB-DEIM algorithm, combining both offline phases into one. Numerical experiments show the capabilities of this novel approach in comparison with classical RB and RB-DEIM approaches.}




\maketitle

\section{Introduction}
\label{sec: intro}

The usual approach for the numerical solution of parabolic partial differential equations (PDEs) is done by time-stepping schemes based upon variational semi-discretizations. After a variational formulation in space and a discretization by, e.g., finite elements, one derives an evolution problem in time. Then, a spatial problem needs to be solved at each discrete time instance.

Alternatively, in a space-time scheme one approximates the PDE by a simultaneous discretization of the spatial \emph{and} temporal domain. This leads to a single variational problem with test functions depending on the spatial and the temporal variables $\bx$ and $t$, respectively.

In this paper, we consider a parameter dependent, nonsmooth parabolic PDE: For a parameter $\mu \in \Pmad$ consider for almost all (f.a.a.) $t \in (0,T]$:
\begin{equation}
\label{eq: PDE}
\tag{$\mu$IVP}
\begin{alignedat}{2}
\dot{y}(t; \mu) - c(\mu) \Delta y(t;\mu) + a(\mu) \max\{0,y(t; \mu)\} & = f(t; \mu) \;&& \text{in } V',\\
y(0; \mu) & = 0 && \text{in } H,
\end{alignedat}
\end{equation}
where $\dot{y} = \frac{\partial y}{\partial t}$ denotes the weak temporal derivative of $y$. Notice that the nonsmoothness in \eqref{eq: PDE} is only Lipschitz continuous and not Fréchet-differentiable. Problem \eqref{eq: PDE} can be seen as a model problem for a broader class of nonsmooth parabolic PDEs, where the assumptions on the nonsmoothness are chosen appropriately (especially to ensure the estimates of Proposition~\ref{pro: stability_FE} and Corollary~\ref{lem: stability_RB}.)

Let us mention some of the related work. Space-time methods have been considered by many authors. We refer, e.g., to the work \cite{Mei11,GK11,NV12,LS15,Ste15,SY18,HT18,Hin20} for (smooth) parabolic problems but there is no error analysis done for nonsmooth PDEs. In the context of reduced basis (RB) methods space-time methods are discussed, e.g., in \cite{SU12,YPU14,Urb14,HPSU22} and in particular for optimal control problems in \cite{BRU22}. A-posteriori error estimates are derived and efficient tensor-based solution methods are proposed. However, the authors do not study empirical interpolation methods and adaptive basis update strategies in their work. Empirical interpolation techniques are necessary to handle the nonsmooth term in \eqref{eq: PDE} efficiently by the RB method, cf. \cite{BMNP04,CS10,CS12}. Theoretical results for \eqref{eq: PDE} can be found in \cite{Bet19,MS17}, where optimal control problems for more general parabolic nonsmooth problems are considered. In the context of adaptive RB methods we refer to \cite{DHO10}, for instance. Let us also mention that the present paper extends results of the elliptic case (cf. \cite{Ber20}) to the parabolic one.

The new contributions of the present paper are as follows: (i) We derive a-posteriori error estimates for space-time approximations of a \emph{nonsmooth} parabolic PDEs. (ii) We incorporate the \emph{discrete emprirical interpolation method} (DEIM) for the nonsmooth term in our error analysis. (iii) A \emph{certified adaptive algorithm} for the RB-DEIM approximation is developed, which combines the typical two offline phases for the computation of an RB and of a DEIM basis into one.

This paper is organized as follows: In Section \ref{sec: space-time}, space-time formulations are introduced for the continuous, finite element (FE) and RB formulations of \eqref{eq: PDE}. Section \ref{sec: RB_method} covers the RB method, including basis generation and space-time a-posteriori error estimation. In Section \ref{sec: DEIM}, DEIM is introduced to efficiently evaluate and approximate the nonsmoothness. The error estimator is adopted to the RB-DEIM setting and an adaptive algorithm for the simultaneous generation of the RB and the DEIM basis is presented. Numerical results illustrate the capabilities this novel approach in Section \ref{sec: numerics}. Finally, we draw some conclusions in Section \ref{sec: conclusion}.
\section{Space-time formulation for \eqref{eq: PDE}}
\label{sec: space-time}

First we establish a space-time formulation for the continuous setting and existence of a unique solution is proved. Then, an FE space-time discretization is introduced. Further, it is shown that its space-time formulation corresponds to a Crank-Nicolson scheme, which can be uniquely solved by a semismooth Newton method. Finally the same is done for the RB space-time formulation. We also present stability estimates for all formulations, that will be necessary to provide convergence rates for some of the RB error estimators quantities.

\subsection{Problem formulation}

Let $\Omega \subset\mathbb R^d$, $d\in\{1,2,3\}$, be a bounded domain with Lipschitz-continuous boundary $\Gamma=\partial\Omega$. We write $\bx=(x_1,\ldots,x_d)$ for an element in $\Omega$. For $T>0$ we define $Q=(0,T)\times\Omega$ and $\Sigma=(0,T)\times\Gamma$. Let $H=L^2(\Omega)$ and $V=H_0^1(\Omega)$ be supplied by the inner products
\begin{align*}
    {\langle\varphi,\phi\rangle}_H&=\int_\Omega\varphi(\bx)\phi(\bx)\, \mathrm d\bx &&\text{for }\varphi,\phi\in H,\\
    {\langle\varphi,\phi\rangle}_V&=\int_\Omega\nabla\varphi(\bx)\cdot \nabla\phi(\bx)\,\mathrm d\bx&&\text{for }\varphi,\phi\in V,
\end{align*}
respectively, and their corresponding induced norms. Moreover, their dual spaces are denoted as $H'$ and $V'$. Furthermore, we introduce the test space $\Ym=L^2(0,T;V)$. We identify the dual $\Ym'$ with the space $L^2(0,T;V')$. For more details on Sobolev and Bochner spaces we refer the reader to \cite{Eva10}, for instance.

Recall that $V \hookrightarrow H \simeq H' \hookrightarrow V'$ is a Gelfand triple and the function space $W(0,T)=\Ym \cap H^1(0,T;V')$ a Hilbert space with induced norm
\begin{align*}
    {\|\varphi\|}_{W(0,T)}^2={\|\varphi\|}_\Ym^2+{\|\dot\varphi\|}_{\Ym'}^2\quad\text{for }\varphi\in W(0,T).
\end{align*}
Moreover, the solution space is $\Xm=\{\varphi \in W(0,T)\,\vert\,\varphi(0) = 0\text{ in }H\}$ with norm 
\begin{align*}
    {\|\varphi\|}_\Xm^2 ={\|\varphi\|}_{W(0,T)}^2 + {\|\varphi(T)\|}_H^2 \quad \text{for }\varphi\in\Xm,
\end{align*}
which is well-defined due to $W(0,T) \hookrightarrow C([0,T]; \Hs)$, cf. \cite{Zei89a}. Finally, we set $\Hm=L^2(0,T;\Hs)$ for brevity. For a function $\varphi\in\Hm$ we write $\varphi(t)$ for the function $\bx\mapsto\varphi(t,\bx)$ for almost all (abbreviated ``f.a.a.'' in the following) $t\in [0,T]$. The precise assumptions on \eqref{eq: PDE} are stated in the following assumption.

\begin{assumption}
    \label{Assumption1}
    \begin{enumerate}
        \item [a)] $\Pmad\subsetneq \mathbb R^p$, $p \in \mathbb N\setminus\{0\}$, is nonempty and compact,
        \item [b)] $c\colon\Pmad\to\mathbb R$ is Lipschitz-continuous, positive and uniformly bounded away from zero,
       \item [c)] $a\colon\Pmad\to\mathbb R$ is Lipschitz-continuous and nonnegative,
       \item [d)] $f\colon \Pmad\to C([0,T];H)$ is Lipschitz-continuous.
    \end{enumerate}
\end{assumption}

\begin{remark}
    Assumption~\ref{Assumption1}-d) can be relaxed to $f\colon\Pmad\to\Hm$. However, we suppose that $t\mapsto f(t;\mu)\in H$ is continuous to simplify our presentation regarding the temporal discretization carried out later.\hfill$\Diamond$
\end{remark}

For $y \in \Xm$, $\phi \in \Ym$ and $\mu \in\Pmad$ we define the operators
\begin{align}
    \label{Definition_A}
    \begin{aligned}
        \calB(y, \phi; \mu)&= \calB_1(y, \phi) + c(\mu) \calB_2(y, \phi),& \calB_1(y,\phi)&={\langle \dot y, \phi\rangle}_{\Ym',\Ym},\\
        \calA(y,\phi;\mu)&= \calB(y,\phi; \mu) + \calN(y, \phi; \mu),
        & \calB_2(y, \phi)&={\langle y,\phi \rangle}_\Ym,\\
        \calN(y, \phi;\mu)&= a(\mu)\,{\langle \max\{0,y\},\phi\rangle}_\Hm,& \calF(\phi;\mu)&={\langle f(\mu),\phi\rangle}_\Hm.
    \end{aligned}
\end{align}
For $\mu \in\Pmad$ the function $y=y(\mu) \in \trial$ is called a \emph{weak solution} to \eqref{eq: PDE} if
\begin{align}
    \label{eq: variational_continuous}
    \calA(y,\phi;\mu) & = \calF(\phi;\mu)\quad \text{for all } \phi \in \Ym.
\end{align}
Existence and uniqueness of the solution to \eqref{eq: variational_continuous} follows e.g. from \cite[Theorem~30.A]{Zei89b}.
\subsection{FE space-time formulation}
\label{sec:FE_space_time}

Analogously to the continuous setting we introduce a discretized space-time formulation. Therefore let $\Xm_\delta\subset\Xm$ and $\Ym_\delta\subset\Ym$ be finite dimensional subspaces. For $\mu \in \Pmad$ we call $y_\delta= y_\delta(\mu) \in \Xm_\delta$ a discretized weak solution to \eqref{eq: PDE} if
\begin{align}
    \label{eq: variational_discretized}
    \calA(y_\delta,\phi_\delta;\mu) & = \calF(\phi_\delta;\mu)\quad \text{for all } \phi_\delta \in\Ym_\delta.
\end{align}
In the remainder of this paper we will focus on the case
\begin{align*}
    \trial_\delta=\FEtimeL \otimes \FEspace,\quad \test_\delta= \FEtimeC \otimes \FEspace,
\end{align*}
where $\otimes$ denotes the tensor product and $\FEtimeC$, $\FEtimeL$ are piecewise constant, respective piecewise linear finite elements in time and $\FEspace$ are piecewise linear finite elements in space. As $\delta=(\Delta t, h)$ we summarize the temporal and spatial discretization parameters. The solutions for this particular choice of spaces will be also called FE solutions to the PDE. For the FE spaces we denote the standard bases (see \cite{Urb14}) by
\begin{align*}
    \FEtimeC = \text{span}\left\{\tau_1, \ldots, \tau_K\right\},\; \FEtimeL = \text{span}\left\{\sigma_1,\ldots,\sigma_K\right\}, \;
    \FEspace  = \text{span}\left\{\zeta_1,\ldots,\zeta_N\right\},
\end{align*}
where $K \in \mathbb N$ and $N\in\mathbb N$ denote the sizes of the temporal and the spatial discretization, respectively. We set $I=[0,T)$ and introduce the intervals $I_k=[t_{k-1}, t_k)$ of length $\Delta t_k$ for the time instances $0 = t_0<\ldots<t_K=T$, i.e., $I = I_1 \cup \ldots \cup I_K$. From
\begin{align*}
    & y_\delta=\sum_{k=1}^K \sum\limits_{i=1}^N \yfec_i^k \big(\sigma_k \otimes \zeta_i\big) \in \Xm_\delta, && \phi_\delta=\sum_{k = 1}^K \sum_{i=1}^N \upphi_i^k \big(\tau_k \otimes \zeta_i\big) \in \test_\delta
\end{align*}
we infer that
\begin{align*}
    \calB(y_\delta,\phi_\delta; \mu) & = \sum_{k,l=1}^K \sum_{i,j = 1}^N \yfec_j^l \upphi_i^k \Big({\langle \dot{\sigma}_l, \tau_k\rangle}_{L^2(0,T)} {\langle \zeta_j, \zeta_i \rangle}_\Hs\\
    &\hspace{28mm} + c(\mu) {\langle \sigma_l, \tau_k \rangle}_{L^2(0,T)} {\langle\zeta_j,\zeta_i \rangle}_V \Big) = \yfec_\delta^\top \Bmat(\mu) \upphi
\end{align*}
with $\yfec_\delta = [\yfec_\delta^1\vert\ldots\vert \yfec_\delta^K]\in\mathbb R^{N\times K}$, $\yfec_\delta^k= (\yfec_1^k, \ldots, \yfec_N^k)^\top$ for $1 \leq k \leq K$, $\upphi=[\upphi^1\vert\ldots\vert\upphi^K]\in\mathbb R^{N\times K}$, $\upphi^k=(\upphi_1^k,\ldots,\upphi_N^k)^\top$ for $1 \leq k \leq K$ and $\Bmat(\mu)=\Bmat_1 + c(\mu) \Bmat_2$. Here we have introduced spatio-temporal matrices $\Bmat_1 =\Ntime \otimes \Mspace$ and $\Bmat_2=\Mtime \otimes \Vspace$ with
\begin{align*}
    & \Mspace &&= \big(\big({\langle\zeta_j, \zeta_i\rangle}_H\big)\big)_{1 \leq i,j\leq N}, && \Vspace &&=\big(\big({\langle\zeta_j,\zeta_i\rangle}_V\big)\big)_{1 \leq i,j \leq N},\\
    & \Mtime &&=\big(\big(\langle \sigma_l, \tau_k\rangle_{L^2(0,T)}\big)\big)_{1 \leq k, l \leq K}, && \Ntime&&=\big(\big({\langle \dot{\sigma}_l, \tau_k \rangle}_{L^2(0,T)}\big)\big)_{1 \leq k,l \leq K}.
\end{align*}
Let $\delta_{l,k}$ denote the Kronecker delta, then we obtain the explicit forms
\begin{align*}
    & {\langle \dot{\sigma}_l, \tau_k \rangle}_{L^2(0,T)} = \delta_{l,k} - \delta_{l+1, k}, && {\langle \sigma_l, \tau_k \rangle}_{L^2(0,T)} = \frac{1}{2} \left(\Delta t_l \delta_{l,k} + \Delta t_{l+1} \delta_{l+1,k}\right).
\end{align*}
Furthermore we denote the lumped version of the spatial mass matrix $\Mspace$ as
\begin{align*}
    \Mlump=\mathrm{diag}\bigg(\frac{1}{3}\big[  \vert\mathrm{supp}(\zeta_j)\vert\big]_{1 \leq j \leq N}\bigg),
\end{align*}
cf., e.g. \cite[Chapter 15]{Tho97}. This leads to
\begin{align*}
    \calB(y_\delta, \tau_k \otimes \zeta_i; \mu) & = \left[\Mspace \left(\yfec_\delta^k - \yfec_\delta^{k-1}\right) + \frac{c(\mu)\Delta t_k}{2}\,\Vspace \left(\yfec_\delta^k + \yfec_\delta^{k-1}\right)\right]_i,\\
    \calN(y_\delta, \tau_k \otimes \zeta_i; \mu) & \approx \left[\frac{a(\mu)\Delta t_k}{2}\,  \Mlump \left(\max\big\{0, \yfec_\delta^k\big\} + \max\big\{0, \yfec_\delta^{k-1}\big\}\right)\right]_i,\\
    \calF(\tau_k \otimes \zeta_i; \mu) & \approx \frac{\Delta t_k}{2} \left(\Fmat_i^k(\mu) + \Fmat_i^{k-1}(\mu)\right)
\end{align*}
for $1\le k\le K$ and $1\le i\le N$, where we have used the trapezoidal quadrature rule for the approximation of the integrals in $\calN$ and $\calF$ and denote $\Fmat(\mu) = [\Fmat^0(\mu)\vert\ldots\vert\Fmat^K(\mu)] \in \mathbb{R}^{N \times K + 1}$ with $\Fmat^k(\mu)=(\Fmat_1^k(\mu),\ldots,\Fmat_{N}^k(\mu))^\top$ and $\Fmat_i^k(\mu) = \langle f(t_k; \mu), \zeta_i\rangle_{\Vs', \Vs}$, where we have utilized that Assumption~\ref{Assumption1}-d) holds. Furthermore we set $\yfec_\delta^0  = 0 \in \mathbb{R}^N$ due to the homogeneous initial condition. Now we can express \eqref{eq: variational_discretized} as a sequence of root finding problems for $k = 1, \ldots, K$:
\begin{align*}
    G_\delta^k(\yfec_\delta^k; \mu)=0\quad\text{in }\mathbb R^N 
\end{align*}
with
\begin{align*}
    &G_\delta^k(\yfec_\delta^k; \mu) =\frac{1}{\Delta t_k} \Mspace \big(\yfec_\delta^k - \yfec_\delta^{k-1}\big) + \frac{1}{2} \big[c(\mu) \Vspace \big(\yfec_\delta^k + \yfec_\delta^{k-1}\big) \\
    &\hspace{10mm} + a(\mu) \Mlump \left(\max\big\{0, \yfec_\delta^k\big\} + \max\big\{0, \yfec_\delta^{k-1}\big\}\right)
- \big(\Fmat^k(\mu) + \Fmat^{k-1}(\mu)\big)\big]
\end{align*}
and initial condition $\yfec_\delta^0=0$. This problem can be interpreted as a Crank-Nicolson (CN) scheme for a spatially discretized parabolic PDE. To solve this problem, a semismooth Newton method (see \cite{Hin10}), where the $k$-th iteration matrix is given by
\begin{align*}
    \mathsf H_\delta^k(\yfec_\delta^k;\mu) =\frac{1}{\Delta t_k} \Mspace + \frac{1}{2} \left[c(\mu) \Vspace + a(\mu) \Mlump \Theta\big(\yfec_\delta^k\big)\right] \in \mathbb R^{N\times N}
\end{align*}
is applied. The function $\Theta \colon \mathbb R^N\to \mathbb R^{N\times N}$ maps a vector to the diagonal matrix that takes the Heaviside function with value $0$ evaluated for each entry of the vector as its diagonal entries.

\begin{remark}
    \label{rem: NewtonSystem}
    For every parameter $\mu \in\Pmad$ the problem $G_\delta^k(\yfec_\delta^k;\mu) = 0$ for $k = 1, \ldots, K$ with initial condition $\yfec_\delta^0 = 0$ admits a unique sequence of roots $\yfec_\delta^k$. This is a consequence of the monotonicity of $G_\delta^k(\cdot\,;\mu)$, which follows since $\Mspace, \Vspace$ and $\Mlump$ are symmetric and positive definite (s.p.d.) matrices and the $\max$-function is monotone. Furthermore for every $\yfec_\delta^k$ the matrix $\mathsf H_\delta^k(\yfec_\delta^k;\mu)$ is s.p.d., which implies that every Newton iteration is uniquely solvable.\hfill$\Diamond$ 
\end{remark}

We finish our considerations with a stability estimate for the FE solution. This is based on the linear case shown in \cite[Corollary 4.3]{Mei11}.

\begin{proposition}
    \label{pro: stability_FE}
    Let Assumption~{\em\ref{Assumption1}} hold. Then, for every $\mu \in \Pmad$, the FE solution $y_\delta(\mu)$ to \eqref{eq: variational_discretized} satisfies the estimate
    \begin{align}
        \label{Estimate1}
        {\|\dot{y}_\delta(\mu)\|}_\Hm\le {\|f(\mu)\|}_\Hm.
    \end{align}
\end{proposition}

\begin{proof}
    Let $y_\delta=y_\delta(\mu)$ be an FE solution to \eqref{eq: variational_discretized}. Since $\dot{y}_\delta\in \Ym_\delta$ holds, we can choose $\phi_\delta = \dot{y}_{\delta}$ as test function. This implies
    \begin{align*}
        {\langle \dot{y}_\delta, \dot{y}_\delta\rangle}_\Hm+ c(\mu)\,{\langle y_\delta,\dot y_\delta \rangle}_\Ym+a(\mu)\,{\langle \max\{0, y_\delta\}, \dot{y}_\delta\rangle}_\Hm= {\langle f(\mu), \dot y_\delta\rangle}_\Hm.
    \end{align*}
    Note that
    \begin{align*}
        c(\mu){\langle y_\delta,\dot y_\delta \rangle}_\Ym&=\frac{c(\mu)}{2} \big({\|y_\delta(T)\|}_V^2 - {\|y_\delta(0)\|}_V^2\big)= \frac{c(\mu)\,{\|y_\delta(T)\|}_V^2}{2} \geq 0,\\
        {\langle f(\mu), \dot{y}_\delta\rangle}_\Hm &\leq\frac{1}{2}\, \big({\|f(\mu)\|}_\Hm^2 + {\|\dot y_\delta\|}_\Hm^2\big).
    \end{align*}
    Thus we obtain
    \begin{align*}
        {\|\dot y_\delta\|}_\Hm^2 + 2 a(\mu)\,{\langle \max\{0, y_\delta\}, \dot y_\delta\rangle}_\Hm \leq {\|f(\mu)\|}_\Hm^2.
    \end{align*}
    Now for almost all $t \in (0,T)$ we introduce the set $\Omega^+(t)= \{\bx \in \Omega\,\vert\,y_\delta(t,\bx) > 0\text{ f.a.a. }t\in[0,T]\}$. Since $\max\{0, \cdot\} \colon \mathbb R \to \mathbb R$ is Lipschitz continuous, we infer from \cite[Chapter~5.8.2.b]{Eva10} that $\max\{0, y_\delta(t)\}$ belongs to $V$ f.a.a. $t\in(0,T)$ and thus
    \begin{align*}
        & \int_I\int_{\Omega^+(t)} y_\delta\dot y_\delta \,\mathrm d\bx\mathrm dt = \int_I \int_\Omega \max\{0, y_\delta\} \dot y_\delta\,\mathrm d\bx\mathrm dt\\
        & \quad = \int_\Omega \max\{0, y_\delta(T)\} y_\delta(T)\,\mathrm d\bx - \int_\Omega \max\{0, y_\delta(0)\} y_\delta(0)\,\mathrm d\bx\\
        & \quad\quad - \int_I \int_\Omega {\max}'\{0, y_\delta\} \dot{y}_\delta y_\delta\,\mathrm d\bx\mathrm dt\\
        & \quad=\int_\Omega \max\{0, y_\delta(T)\} y_\delta(T)\,\mathrm d\bx - \int_I \int_{\Omega^+(t)} y_\delta \dot{y}_\delta\,\mathrm d\bx\mathrm dt,
    \end{align*}
    since for the weak derivative $\max'\{0, y_\delta\} = 1$ almost everywhere (abbreviated ``a.e.'' in the following) on $\Omega^+(t)$ and $\max'\{0, y_\delta\} = 0$ a.e. on $\Omega \setminus \Omega^+(t)$ are satisfied; cf. \cite{Ber19} for a derivation based on generalized derivatives. This implies
    \begin{align*}
        {\langle \max\{0, y_\delta\}, \dot{y}_\delta\rangle}_\Hm = \frac{1}{2} {\langle \max\{0, y_\delta(T)\}, y_\delta(T)\rangle}_H \geq 0
    \end{align*}
    and thus \eqref{Estimate1} holds true.
\end{proof}

For an arbitrarily given $\varphi_\delta\in\Xm_\delta$, which is piecewise linear in time, let us define the piecewise constant in time version by averaging
\begin{align}
    \label{phik_all}
    \bar\varphi_\delta(t,\bx)= \sum_{k=1}^K \mathbbm{1}_{I_k}(t)\otimes \bar\varphi_\delta^k(\bx)\in\Ym_\delta
\end{align}
with
\begin{align}
    \label{phik}
    \bar\varphi_\delta^k= \frac{1}{\Delta t_k} \int_{I_k} \varphi_\delta(t)\, \mathrm dt=\frac{\varphi(t_{k-1})+\varphi(t_k)}{2}\in V. 
\end{align}

In the next lemma we summarize useful properties of $\dot\varphi_\delta \in \Hm$. Recall that $\dot\varphi_\delta \in \Hm$ for $\varphi_\delta\in \Xm_\delta$.

\begin{lemma}
    \label{lem:W0T_property_disc_and_ProjOrtho}
    For every $\varphi_\delta\in\Xm_\delta$ we have
    \begin{align}
        \label{Lemma:Est1}
        {\langle \dot\varphi_\delta, \bar\varphi_\delta\rangle}_\Hm= \frac{1}{2} {\|\varphi_\delta(T)\|}_H^2,
    \end{align}
    and
    \begin{align}
        \label{Lemma:Est2}
        {\langle\varphi_\delta-\bar\varphi_\delta,\phi_\delta\rangle}_\Ym={\langle\varphi_\delta-\bar\varphi_\delta,\phi_\delta\rangle_\Hm}=0\quad\text{for all }\phi_\delta\in\Ym_\delta.
    \end{align}
\end{lemma}

\begin{proof}
    Let $\varphi_\delta\in\Xm_\delta$ be chosen arbitrarily. Equation \eqref{Lemma:Est1} follows from
    \begin{align*}
        {\langle \dot\varphi_\delta, \bar\varphi_\delta\rangle}_\Hm &= \sum_{k = 1}^K \int_{I_k} \bigg\langle \frac{\varphi_\delta(t_k)-\varphi_\delta(t_{k-1})}{\Delta t_k}, \frac{\varphi_\delta(t_k)+\varphi_\delta(t_{k-1})}{2} \bigg\rangle_H\,\mathrm dt\\
        & = \frac{1}{2}\sum_{k=1}^K \big({\|\varphi_\delta(t_k)\|}_H^2 - {\|\varphi_\delta(t_{k-1})\|}_H^2\big)=\frac{1}{2} {\|\varphi_\delta(T)\|}_H^2,
    \end{align*}
    where we have used that $\varphi_\delta(0)=0$ in $H$.\hfill\\
    To show \eqref{Lemma:Est2} let $\varphi_\delta\in\Xm_\delta$ and $\phi_\delta\in\Ym_\delta$ be chosen arbitrarily. Then, it follows that
    \begin{align*}
        \varphi_\delta(t)&=\frac{1}{\Delta t_k} \,\big((t_k - t)\varphi_\delta(t_{k-1})+ (t - t_{k-1})\varphi_\delta(t_k)\big)&&\text{f.a.a. }t\in I_k,\\
        \bar\varphi_\delta(t)&=\frac{\varphi_\delta(t_{k-1})+\varphi_\delta(t_k)}{2}&&\text{f.a.a. }t\in I_k,\\
        \phi_\delta(t)& = \phi_\delta(t_k) =:\phi_\delta^k&&\text{f.a.a. }t\in I_k
    \end{align*}
    for every $k=1,\ldots,K$. Consequently,
    \begin{align*}
        {\langle\varphi_\delta,\phi_\delta\rangle}_\Ym&= \frac{1}{\Delta t_k} \sum_{k=1}^K \bigg\langle\Big[\Big(t_kt - \frac{t^2}{2}\Big)\varphi_\delta(t_{k-1})+\Big(\frac{t^2}{2} - t_{k-1} t\Big)\varphi_\delta(t_k),\phi_\delta^k\Big]_{t=t_{k-1}}^{t=t_k} \bigg\rangle_V\\
        &= \sum_{k=1}^K \frac{1}{\Delta t_k} \bigg\langle \frac{(t_k - t_{k-1})^2}{2} \big((\varphi_\delta(t_{k-1})+\varphi(t_k)\big),\phi_\delta^k\bigg\rangle_V\\
        &= \sum\limits_{k = 1}^K \frac{\Delta t_k}{2} {\langle\varphi_\delta(t_{k-1})+\varphi_\delta(t_k),\phi_\delta^k\rangle}_V
    \end{align*}
    and
    \begin{align*}
        {\langle\bar\varphi_\delta,\phi_\delta\rangle}_\Ym&=\sum_{k=1}^K\int_{I_k}\bigg\langle \frac{\varphi_\delta(t_{k-1})+\varphi_\delta(t_k)}{2},\phi_\delta^k\bigg\rangle_V\,\mathrm dt\\
        &= \sum_{k = 1}^K \frac{\Delta t_k}{2} {\langle\varphi_\delta(t_{k-1})+\varphi_\delta(t_k),\phi_\delta^k\rangle}_V,
    \end{align*}
    which implies $\langle\varphi_\delta-\bar\varphi_\delta,\phi_\delta\rangle_\Ym=0$. Analogously, we derive that $\langle\varphi_\delta-\bar\varphi_\delta,\phi_\delta\rangle_\Hm=0$. This finishes the proof.
\end{proof}

\subsection{RB space-time formulation}
\label{sec: space-time_RB}

Analogously the space-time RB setting can be formulated and interpreted as a CN scheme. This CN interpretation guides us in the numerical calculation of solutions to the space-time formulation and in the derivation of an a-posteriori error estimator. For a spatial RB space $\RBspace= \text{span}\{\psi_1,\ldots,\psi_\ell\} \subset \FEspace$ of dimension $\ell\in\mathbb N$, we introduce the RB solution and test spaces
\begin{align*}
    & \Xm_\RB=\FEtimeL \otimes \RBspace, &&\test_\RB= \FEtimeC \otimes \RBspace,
\end{align*}
respectively, where $\RB=(\Delta t,\ell)$ stands for the temporal discretization and RB parameter.

For $\mu \in\Pmad$ we call $y_\RB=y_\RB(\mu) \in \trial_\RB$ an RB solution to \eqref{eq: variational_discretized} if
\begin{align}
    \label{eq: variational_RB}
    & \calA(y_\RB,\phi;\mu) = \calF(\phi;\mu) && \text{for all }\phi \in \test_\RB.
\end{align}
As previously done in the case of the FE space-time formulation, we can also reformulate \eqref{eq: variational_RB} as a CN scheme, by using a trapezoidal quadrature rule. Therefore let $\RBmat =[\uppsi_1, \ldots, \uppsi_\ell] \in \mathbb R^{N\times\ell}$ denote the RB coefficient matrix whose columns are the FE coefficient vectors of the reduced basis functions. Let $\MspaceL=\RBmat^T \Mspace \RBmat \in \mathbb R^{\ell \times \ell}$ and $\VspaceL=\RBmatT \Vspace \RBmat \in \mathbb R^{\ell \times \ell}$ denote the spatial RB mass matrix and spatial RB stiffness matrix, respectively. This leads to a sequence of root finding problems for $k = 1,\ldots, K$:
\begin{align*}
    G_\RB^k(\yfec_\RB^k; \mu)=0\quad\text{in }\mathbb R^\ell
\end{align*}
with
\begin{align*}
    G_\RB^k(\yfec_\RB^k; \mu)&= \frac{1}{\Delta t_k} \MspaceL \big(\yfec_\RB^k - \yfec_\RB^{k-1}\big) + \frac{1}{2} \big[ c(\mu) \VspaceL \big( \yfec_\RB^k + \yfec_\RB^{k-1}\big)\\
    &\quad + a(\mu) \RBmatT \Mlump \left(\max\big\{0, \RBmat \yfec_\RB^k\big\} + \max\big\{0, \RBmat \yfec_\RB^{k-1}\big\}\right)\\
    &\quad  - \RBmatT \big(\Fmat^k(\mu) + \Fmat^{k-1}(\mu)\big)\big]
\end{align*}
with initial condition $\yfec_\RB^0 = 0$. Again this problem can be solved by applying a semismooth Newton method, where the $k$-th iteration matrix is given by
\begin{align*}
    \mathsf H_\RB^k(\yfec_\RB^k; \mu)&=\frac{1}{\Delta t_k} \MspaceL\\
    &\quad+\frac{1}{2} \Big(c(\mu) \VspaceL+ a(\mu) \RBmatT \Mlump \Theta\big(\RBmat \yfec_\RB^k\big) \RBmat\Big) \in \mathbb R^{\ell \times \ell}.
\end{align*}
Since $\RBmat$ has full rank, Remark \ref{rem: NewtonSystem} is applicable again, i.e. the sequence of roots $\yfec_\RB^k$ is unique and every Newton iteration is uniquely solvable. Furthermore we obtain a stability estimate analogous to Proposition~\ref{pro: stability_FE}.

\begin{corollary}
    \label{lem: stability_RB}
    Let Assumption~{\em\ref{Assumption1}} hold. Then, for every $\mu\in\Pmad$, the RB solution $y_\RB(\mu)$ to \eqref{eq: variational_RB} satisfies the estimate
    \begin{align*}
        {\|\dot y_\RB(\mu)\|}_\Hm \leq\,{\|f(\mu)\|}_\Hm.
    \end{align*}
\end{corollary}

\begin{proof}
    The claim follows by similar arguments utilized to prove of Proposition~\ref{pro: stability_FE}.
\end{proof}
\section{Reduced basis method}
\label{sec: RB_method}

In this section a greedy procedure for the generation of the RB space $\RBspace$ will be presented. Efficient error estimation is necessary and an error estimator is derived based on the space-time formulations introduced in Section \ref{sec: space-time}. Compared to the elliptic case, cf. \cite{Ber20}, the derived error estimator is composed of an additional term $\Delta_{\calP^\delta}$. A large part of this section is dedicated to its convergence analysis, ultimately proving that $\Delta_{\calP^\delta} = \mathcal{O}(\Delta t)$, where an equidistant temporal discretization is assumed for ease of presentation.

\subsection{Generation of reduced basis}
\label{sec: RB_basis}

For the generation of the spatial RB space $\RBspace$ proper orthogonal decomposition (POD) is used in a greedy procedure. For more details on POD we refer, e.g., to \cite{Vol17}. With $\calP^\ell:V\to V_\ell$ the $V$-orthogonal projection onto the spatial RB space is denoted. By $\Delta(\mu)$ we denote an error estimator for the RB solution corresponding to the parameter $\mu \in \mathscr P$ with respect to RB solution space $\Xm_\RB$ and RB test space $\Ym_\RB$ generated by $V_\ell$. Furthermore $\POD_1$ denotes the extraction of a dominant POD mode with respect to the inner product in $\Vs$. The offline basis generation is summarized in Algorithm \ref{algo: POD_Greedy_RB}.

\begin{algorithm}[h!]
    \caption{(POD-greedy RB method)}\label{algo: POD_Greedy_RB}
    \begin{algorithmic}[1]
        \Require{Discrete training set $\train \subset\Pmad$, error tolerance $\Totaltol > 0$;}
        \State{Set $\ell = 0$, $\mathsf{\Psi}_0 =  \emptyset$, $\Vs_{0} =  \{0\}$;}
        \While{$\varepsilon_\ell:=\max\{\Delta(\mu)\,\vert\,\mu \in \train\} > \Totaltol$}
            \State{Compute $\mu_{\ell+1} \in \argmax\{\Delta(\mu)\,\vert\,\mu \in \train\}$;}
			\State{Set $e_{\ell + 1}^k = y_\delta^k(\mu_{\ell + 1}) - \calP^\ell y_\delta^k(\mu_{\ell + 1})$ for $k = 1,\ldots,K$;}
            \State{Define $\mathscr P_{\ell + 1} =\mathscr P_\ell\cup\{\mu_{\ell + 1}\}$;}
			\State{Compute $\psi_{\ell + 1} \in \POD_1(\{e_{\ell + 1}^k\}_{k = 1}^K)$;}
			\State{Set $\mathsf{\Psi}_{\ell + 1} = \RBmat \cup \{\psi_{\ell + 1}\}$, $V_{\ell + 1} = \RBspace \oplus \text{span}(\psi_{\ell + 1})$ and $\ell = \ell + 1$;}
        \EndWhile\hfill\\
        \Return{Reduced basis $\RBmat$, spatial RB space $\Vs_\ell$.}
    \end{algorithmic}
\end{algorithm}


It is possible to add more than just one POD mode per iteration to the reduced basis. This usually results in faster convergence of the algorithm and thus lower offline computational cost. But on the other hand it does not guarantee that the reduced basis is of minimal size. Note that the initial condition, which is zero, is always perfectly approximated.

\subsection{RB error estimation}

An efficient error estimator is crucial not only for Algorithm \ref{algo: POD_Greedy_RB}, but also for quality certification of the calculated RB solution in an online phase. Since $y_\delta, y_\RB \notin \test_\RB$ we cannot test with the error $y_\delta - y_\RB$ in the derivation of an error estimator. Therefore we will use an additional projection operator $\calP^\delta$ that satisfies the following hypotheses.

\begin{assumption}
    \label{AssumptionProj}
    Let $\calP^\delta\colon \trial_\delta \rightarrow \test_\delta$ be a bounded projection operator such that:
    \begin{enumerate}
        \item [a)] $\calP^\delta$ satisfies
        \begin{align*}
            {\langle\dot \varphi_\delta,\calP^\delta\varphi_\delta\rangle}_\Hm=\frac{1}{2}\,{\|\varphi_\delta(T)\|}_H^2\quad\text{for all }\varphi_\delta\in \trial_\delta.
        \end{align*}
        \item [b)] $\calP^\delta$ satisfies
        \begin{align*}
            {\langle \varphi_\delta-\calP^\delta\varphi_\delta,\phi_\delta\rangle}_\Ym={\langle \varphi_\delta-\calP^\delta\varphi_\delta,\phi_\delta\rangle}_\Hm=0\quad\text{for all }\varphi_\delta\in\Xm_\delta,\phi_\delta\in \Ym_\delta.
        \end{align*}
    \end{enumerate}
\end{assumption}

We have already seen an example of a suitable projection operator in the previous section.

\begin{example}
    For $1 \leq k \leq K$ and $t \in I_k$ almost everywhere we introduce the projection
    \begin{align*}
        \big(\calP^\delta \varphi_\delta\big)(t) =\bar\varphi_\delta^k \quad \text{for all } \varphi_\delta\in \Xm_\delta,
    \end{align*}
    where $\bar\varphi_\delta^k$ is defined in \eqref{phik_all}. Due to Lemma~\ref{lem:W0T_property_disc_and_ProjOrtho},  Assumption~\ref{AssumptionProj} is satisfied.\hfill$\Diamond$
\end{example}
 
 Before stating the main theorem of this section, we summarize properties of the projection operator $\calP^\delta$ in the next lemma.
 
\begin{lemma}
    \label{lem: FE_projection1}
    Suppose that the projection operator $\calP^\delta$ satisfies Assumption~\em{\ref{AssumptionProj}-b)}. For $\varphi_\delta\in\Xm_\delta$ it holds that
    \begin{enumerate}
        \item [a)] $\|\varphi_\delta- \calP^\delta\varphi_\delta\|_\Hm \leq \|\varphi_\delta-\phi_\delta\|_\Hm$ for all $\phi_\delta\in \Ym_\delta$.
        \item [b)] $\|\varphi_\delta-\calP^\delta \varphi_\delta\|_\Hm\le\Delta t\,\|\dot\varphi_\delta\|_\Hm$, where $\Delta t=\max_{1\le k\le K}\Delta t_k$ denotes the maximal time step.
    \end{enumerate}
\end{lemma}

\begin{proof}
    \begin{enumerate}
        \item [a)] For $\varphi_\delta\in \Xm_\delta$ and $\phi_\delta\in \Ym_\delta$ Assumption~\ref{AssumptionProj}-2) and $\calP^\delta\varphi_\delta\in\Ym_\delta$ imply that
        \begin{align*}
            {\|\varphi_\delta-\calP^\delta \varphi_\delta\|}_\Hm^2& = {\langle \varphi_\delta-\calP^\delta\varphi_\delta,\varphi_\delta\rangle}_\Hm-{\langle \varphi_\delta-\calP^\delta\varphi_\delta,\calP^\delta\varphi_\delta\rangle}_\Hm\\
            &= {\langle \varphi_\delta-\calP^\delta\varphi_\delta,\varphi_\delta\rangle}_\Hm\\
            & = {\langle \varphi_\delta- \calP^\delta\varphi_\delta,\varphi_\delta-\phi_\delta \rangle}_\Hm\le{\|\varphi_\delta- \calP^\delta\varphi_\delta\|}_\Hm{\|\varphi_\delta-\phi_\delta\|}_\Hm,
        \end{align*}
        which gives the claim.
        \item [b)] Denote by $\Pi^\delta \colon \trial_\delta \rightarrow \test_\delta$ the piecewise constant interpolation operator defined as $(\Pi^\delta \varphi_\delta)(t)= \varphi_\delta(t_{k-1})$ for $1 \leq k \leq K$ such that $t \in I_k$. For $e(t)=\varphi_\delta(t)-(\Pi^\delta \varphi_\delta)(t) \in \FEspace$, we denote by $e_k$ the restriction of $e$ on the interval $I_k$. Then $e_k$ is affine linear w.r.t. $t$ and belongs to $H^1(I_k;V)$. Thus, for $t \in I_k$ we estimate
        \begin{align*}
            e(t)=e_k(t) = \underbrace{e_k(t_{k-1})}_{=\, 0}+\int_{t_{k-1}}^t \dot{e}_k(s)\,\mathrm{d}s
        \end{align*}
        by the fundamental theorem of calculus. This implies
        \begin{align*}
            {\|e_k(t)\|}_H & = \bigg\| \int_{t_{k-1}}^t \dot{e}_k(s)\,\mathrm ds \bigg\|_H\leq\int_{I_k} {\|\dot{e}_k(s)\|}_H \,\mathrm ds\\
            &\leq\bigg(\int_{I_k} 1 \,\mathrm ds\bigg)^{1/2} \bigg(\int_{I_k} {\|\dot e_k(s)\|}_H^2\,\mathrm ds \bigg)^{1/2}=\sqrt{\Delta t_k}\,{\|\dot e_k\|}_{L^2(I_k; H)}.
        \end{align*}
        Consequently,
        \begin{align*}
            {\| e_k\|}_{L^2(I_k;H)}^2 & = \int_{I_k} {\| e_k(s)\|}_H^2\,\mathrm ds\leq \left(\Delta t_k\right)^2 {\| \dot e_k \|}_{L^2(I_k;H)}^2.
        \end{align*}
        Since $\dot e_k= \dot\varphi_\delta$ on $I_k$ this implies that
        \begin{align*}
            {\|e_k\|}_{L^2(I_k;H)} \leq \Delta t_k \,{\|\dot \varphi_\delta\|}_{L^2(I_k;H)}.
        \end{align*}
        We have $\Pi^\delta\varphi_\delta\in\Ym_\delta$ and $\Delta t=\max_{1\le k\le K}\Delta t_k$. Hence, using part a) it follows
        \begin{align*}
            {\|\varphi_\delta-\calP^\delta \varphi_\delta\|}_\Hm^2&\leq {\|\varphi_\delta-\Pi^\delta \varphi_\delta\|}_\Hm^2 = \sum_{k = 1}^K {\|\varphi_\delta - \Pi^\delta \varphi_\delta\|}_{L^2(I_k;H)}^2\\
            &\leq\sum_{k=1}^K (\Delta t_k)^2 {\|\dot\varphi_\delta\|}_{L^2(I_k;H)}^2\le(\Delta t)^2\,{\|\dot\varphi_\delta\|}_\Hm^2.
        \end{align*}
        This implies part b).
    \end{enumerate}
\end{proof}

Next we define the RB residual as
\begin{align}
    \label{Residuum}
    \calR_\RB(\phi;\mu)=\calF(\phi;\mu)-\calA(y_\RB(\mu),\phi;\mu) \quad \text{for all }\phi \in \test_\delta.
\end{align}
Now we are ready to state the RB error estimator.

\begin{theorem}
    \label{thm: error_estimator_fix_new}
    Let Assumption~{\em\ref{Assumption1}} hold. Suppose that the projection operator $\calP^\delta$ satisfies Assumption~{\em\ref{AssumptionProj}}. For $\mu \in\Pmad$ denote the RB error as $e_\RB(\mu)=y_\delta(\mu)-y_\RB(\mu)$ satisfying $\calP^\delta e_\RB(\mu)\neq0$. Then,
    \begin{align}
        \label{Estimate2}
        {\|\calP^\delta e_\RB(\mu)\|}_\Ym & \leq \Delta_\RB(\mu)+\Delta_{\calP^\delta}(\mu),
    \end{align}
    with RB error estimator $\Delta_\RB$ and projection error estimator $\Delta_{\calP^\delta}$ given by
    \begin{align*}
        \Delta_\RB(\mu)&= \frac{1}{c(\mu)}\,{\| \calR_\RB(\cdot\,;\mu)\|}_{\Ym_\delta'},\\
        \Delta_{\calP^\delta}(\mu) &=\frac{a(\mu)}{c(\mu)} \,\frac{{\|\max\{0, y_\delta(\mu)\} - \max\{0, y_\RB(\mu)\}\|}_\Hm}{{\|\calP^\delta e_\RB(\mu)\|}_\Ym}\,{\|e_\RB(\mu) - \calP^\delta e_\RB(\mu)\|}_\Hm,
    \end{align*}
    respectively.
\end{theorem}

\begin{proof}
    Let $\mu\in\Pmad$ be chosen arbitrarily. From $e_\RB(\mu)\in\Xm_\delta$ we obtain by Assumption~\ref{AssumptionProj}-b), the monotonicity of $\max$ and Assumption~\ref{AssumptionProj}-a) that
    \begin{align*}
        &{\|\calP^\delta e_\RB(\mu)\|}_\Ym^2 ={\langle\calP^\delta e_\RB(\mu),\calP^\delta e_\RB(\mu)\rangle}_\Ym={\langle e_\RB(\mu),\calP^\delta e_\RB(\mu)\rangle}_\Ym,\\
        & a(\mu)\,{\langle \max\{0, y_\delta(\mu)\}-\max\{0, y_\RB(\mu)\}, e_\RB(\mu)\rangle}_\Hm \ge 0,\\
        & {\langle \dot e_\RB(\mu), \calP^\delta e_\RB(\mu)\rangle}_\Hm= {\langle \dot{e}_\RB(\mu), e_\RB(\mu)\rangle}_\Hm= \frac{1}{2} {\|e_\RB(T;\mu)\|}_H^2 \geq 0.
    \end{align*}
    Thus we infer that
    \begin{align*}
        & c(\mu)\,{\|\calP^\delta e_\RB(\mu)\|}_\Ym^2 =c(\mu)\,{\langle e_\RB(\mu),\calP^\delta e_\RB(\mu)\rangle}_\Ym\\
        & \leq c(\mu)\,{\langle e_\RB(\mu), \calP^\delta e_\RB(\mu) \rangle}_\Ym+ a(\mu){\langle\max\{0,y_\delta(\mu)\} -\max\{0, y_\RB(\mu)\},e_\RB(\mu) \rangle}_\Hm\\
        & \leq c(\mu)\,{\langle e_\RB(\mu),\calP^\delta e_\RB(\mu)\rangle}_\Ym+{\langle\dot e_\RB(\mu),\calP^\delta e_\RB(\mu)\rangle}_\Hm\\ 
        & \quad+a(\mu)\big({\langle\max\{0, y_\delta(\mu)\}-\max\{0, y_\RB(\mu)\}, \calP^\delta e_\RB(\mu)\rangle}_\Hm\\
        & \quad\hspace{12mm}+{\langle \max\{0, y_\delta(\mu)\}-\max\{0, y_\RB(\mu)\},e_\RB(\mu) - \calP^\delta e_\RB(\mu)\rangle}_\Hm\big).
    \end{align*}
    Utilizing \eqref{Definition_A} it follows that
    \begin{align*}
        c(\mu)\,{\|\calP^\delta e_\RB(\mu)\|}_\Ym^2 & \le a(\mu){\langle \max\{0, y_\delta(\mu)\} - \max\{0, y_\RB(\mu)\}, e_\RB(\mu) - \calP^\delta e_\RB(\mu)\rangle}_\Hm\\
        & \quad+\mathcal A(y_\delta(\mu),\calP^\delta e_\RB(\mu);\mu)-\mathcal A(y_\RB(\mu),\calP^\delta e_\RB(\mu);\mu).
    \end{align*}
    Moreover, we deduce from \eqref{eq: variational_discretized} that
    \begin{align*}
        \mathcal A(y_\delta(\mu),\calP^\delta e_\RB(\mu);\mu)=\calF\big(\calP^\delta e_\RB(\mu);\mu\big).
    \end{align*}
    Consequently, by using \eqref{Residuum}
    \begin{align*}
        & c(\mu)\,{\|\calP^\delta e_\RB(\mu)\|}_\Ym^2 \leq \calR_\RB\big(\calP^\delta e_\RB(\mu);\mu\big)\\
        & \quad+a(\mu){\langle \max\{0,y_\delta(\mu)\}-\max\{0, y_\RB(\mu)\}, e_\RB(\mu)-\calP^\delta e_\RB(\mu)\rangle}_\Hm\\
        & \leq {\| \calR_\RB(\cdot\,;\mu)\|}_{\test_\delta'} {\|\calP^\delta e_\RB(\mu)\|}_\test \\
        & \quad + a(\mu) {\|\max\{0, y_\delta(\mu)\} - \max\{0,y_\RB(\mu)\}\|}_\Hm {\|e_\RB(\mu) - \calP^\delta e_\RB(\mu)\|}_\Hm.
    \end{align*}
    Thus we obtain
    \begin{align*}
        {\|\calP^\delta e_\RB(\mu)\|}_\Ym & \leq \frac{1}{c(\mu)}{\| \calR_\RB(\cdot\,;\mu)\|}_{\Ym_\delta'}\\
        & + \frac{a(\mu)}{c(\mu)} \frac{{\|\max\{0,y_\delta(\mu)\} - \max\{0, y_\RB(\mu)\}\|}_\Hm}{{\|\calP^\delta e_\RB(\mu)\|}_\Ym} {\|e_\RB(\mu) - \calP^\delta e_\RB(\mu)\|}_\Hm
    \end{align*}
    which gives the claim.
\end{proof}

Finally we study the convergence of the projection error estimator $\Delta_{\calP^\delta}(\mu)$.

\begin{proposition}
    \label{prop: time_convergence}
    Let Assumptions~{\em\ref{Assumption1}} and {\em\ref{AssumptionProj}} hold. Denote $\Delta t = \max_{1 \leq k \leq K} \Delta t_k$. Then $\Delta_{\calP^\delta}(\mu) \to 0$ for $\Delta t \to 0$ and every $\mu \in\Pmad$. Especially for sufficiently small $\Delta t$ and every $\mu \in \Pmad$ the estimate
    \begin{align*}
        \Delta_{\calP^\delta}(\mu) \leq C \Delta t
    \end{align*}
    with a constant $C \geq 0$ independent of $\mu$ is satisfied.
\end{proposition}

\begin{proof}
    Without loss of generality, we can assume that $e_\RB(\mu) \neq 0$ for the parameters $\mu$ under consideration, or there is nothing to show. Recall that $e_\RB(\mu)=y_\delta(\mu)-y_\RB(\mu)\in\Xm_\delta$ holds true. By Lemma~\ref{lem: FE_projection1}-b), Proposition \ref{pro: stability_FE} and Corollary~\ref{lem: stability_RB} we obtain
    \begin{align}
    \label{Pconv1}
        \begin{aligned}
            {\|e_\RB(\mu)-\calP^\delta e_\RB(\mu)\|}_\Hm& \le{\|y_\delta(\mu)-\calP^\delta y_\delta(\mu)\|}_\Hm+{\| y_\RB(\mu)-\calP^\delta y_\RB(\mu)\|}_\Hm\\
            & \leq\big({\|\dot y_\delta(\mu)\|}_\Hm+{\|\dot y_\RB(\mu)\|}_\Hm\big) \Delta t\leq 2\,{\|f(\mu)\|}_\Hm \Delta t
        \end{aligned}
    \end{align}
    for every $\mu\in\Pmad$. Thus $\calP^\delta e_\RB(\mu) \to e_\RB(\mu)$ in $\Hm$ for $\Delta t \to 0$. This implies that there exists a (sufficiently small) constant $\tau_\mu\in(0,T]$ satisfying
    \begin{align}
        \label{Pconv2}
        \begin{aligned}
            0\le\frac{{\|e_\RB(\mu)\|}_\Hm}{{\|\calP^\delta e_\RB(\mu)\|}_\Hm}\le\frac{3}{2}\quad\text{for }\Delta t\in\big(0,\tau_\mu\big].
        \end{aligned}
    \end{align}
    Here we have used that without loss of generality $\calP^\delta e_{\RB}(\mu) \neq 0$ for all $\Delta t \in \big(0,\tau_\mu\big]$, possibly after shrinking $\tau_\mu$, because $\calP^\delta e_\RB(\mu) \to e_\RB(\mu)$ in $\Hm$ for $\Delta t \to 0$ and $e_\RB(\mu) \neq 0$ holds. From Poincaré's inequality we infer that there exists a constant $c_P>0$ so that
    \begin{align}
        \label{Pconv3}
        \begin{aligned}
            \frac{1}{{\|\varphi\|}_\Ym}\le \frac{c_P}{{\|\varphi\|}_\Hm}\quad\text{for all }\varphi\in\Ym\setminus\{0\}.
        \end{aligned}
    \end{align}
    Utilizing \eqref{Pconv1}-\eqref{Pconv3} we obtain for $\Delta t\in(0,\tilde C(\mu)]$
    \begin{align*}
        \Delta_{\calP^\delta}(\mu) & = \frac{ a(\mu)}{c(\mu)} \frac{{\|e_\RB(\mu)\|}_\Hm}{{\|\calP^\delta e_\RB(\mu)\|}_\Ym} {\|e_\RB(\mu) - \calP^\delta e_\RB(\mu)\|}_\Hm\\
        &\leq\frac{2a(\mu)}{c(\mu)} \frac{{\|e_\RB(\mu)\|}_\Hm}{{\|\calP^\delta e_\RB(\mu)\|}_\Ym}{\|f(\mu)\|}_\Hm \Delta t \leq\hat C(\mu)\, \Delta t,
    \end{align*}
    where the nonnegative constant $\hat C(\mu)=3 a(\mu)c_P\,\|f(\mu)\|_\Hm/c(\mu)$ is uniformly bounded w.r.t. $\mu$ due to the Lipschitz continuity of $a$, $c$, $f$ and the compactness of $\Pmad$.
\end{proof}
\section{Adaptive RB-DEIM method}
\label{sec: DEIM}

The current model order reduction approach is unsatisfactory for two reasons. First, the sequence of root finding problems $G_\RB^k(\yfec_\RB^k;\mu)$ is on spatial RB level $\ell$, but an evaluation of the nonsmoothness on spatial FE level $N$ is necessary. Furthermore the dual norm of the residual for the RB error estimator $\Delta_\RB(\mu)$ cannot be efficiently computed in an offline-online separable fashion. This is typical for RB error estimation in the context of nonlinear PDEs, cf. \cite{Ber20, Hin20}. We will use the discrete empirical interpolation method (DEIM) to approximate the nonsmoothness and overcome those difficulties. For details on the efficient evaluation of the dual norm of the residual in combination with (D)EIM, we refer to \cite[Section~4.2.5]{Hesthaven}.

In the first half of this section the classical DEIM framework is presented together with a new approach for error estimation. Afterwards an adaptive DEIM framework that combines RB and DEIM offline phases is presented. 

\subsection{Classical DEIM}
\label{sec: ClassicalDEIM}

We keep the presentation of the classical DEIM approach short and refer to \cite{CS10, CS12} for more details on the general procedure and to \cite{Ber19} for an application to an elliptic max-type PDE.
For the remainder of this section we will assume that the right-hand sides $\mathsf{F}^k$ are parameter separable independently of the time instance, i.e., 
\begin{align*}
    \Fmat^k(\mu) = \Bmat_\Fmat \beta_\Fmat(\mu) \gamma_\Fmat(k),
\end{align*}
with
\begin{align*}
    \Bmat_\Fmat \in \mathbb R^{N\times p}, \quad \beta_\Fmat \colon \mathscr P \to \mathbb R^p, \quad \gamma_\Fmat \colon\{0,\ldots,K\} \to \mathbb R.
\end{align*}
If this should not be the case, a further DEIM approximation would be necessary. The procedure can be easily generalized to this situation, cf. \cite{CS10}. Now DEIM introduces a projection matrix $\Pmat = [e_{i_1}\vert\ldots\vert e_{i_L}] \in \mathbb R^{N \times L}$ and an approximation matrix $\Phimat = [\phi_1\vert\ldots\vert\phi_L] \in \mathbb R^{N \times L}$. We assume $L \ll N$ and that $\Pmat^\top \Phimat \in \mathbb R^{L \times L}$ is invertible. Then the max term in $G_\RB^k$ can be approximated as
\begin{align*}
    \Phimat (\Pmat^\top \Phimat)^{-1} \max\big\{0, \Pmat^\top \RBmat \yfec_\RB^k(\mu)\big\}.
\end{align*} 
This gives the RB-DEIM approximation of the RB root finding problem presented in Section \ref{sec: space-time_RB}. For $k = 1,\ldots, K$ solve the sequence of root finding problems:
\begin{align}
    \label{RootRB}
    G_{\RB, L}^k(\yfec_{\RB,L}^k;\mu)=0\quad\text{in }\mathbb{R}^\ell
\end{align}
with
\begin{align*}
     G_{\RB, L}^k(\yfec_{\RB,L}^k;\mu) &=\frac{1}{\Delta t_k} \MspaceL (\yfec_{\RB, L}^k - \yfec_{\RB,L}^{k - 1})+ \frac{1}{2}\Big[c(\mu)\,\VspaceL (y_{\RB,L}^k + y_{\RB,L}^{k - 1})\\
    &\hspace{18mm}+ a(\mu)\, \mathsf\Psi^\top_\ell\Mlump \Phimat (\Pmat^\top \Phimat)^{-1} \big(\max\big\{0, \Pmat^\top \RBmat \yfec_{\RB,L}^k\big\}\\
    &\hspace{18mm}+ \max\big\{0, \Pmat^\top \RBmat \yfec_{\RB,L}^{k - 1}\big\}\big) -\mathsf\Psi^\top_\ell (\mathsf{F}^k(\mu) + \mathsf{F}^{k-1}(\mu))\Big],
\end{align*}
with initial condition $\yfec_{\RB,L}^0 = 0$. Again this problem can be solved by applying a semismooth Newton method, where the $k$-th iteration matrix is given by
\begin{align*}
    &\mathsf H_{\RB,L}^k(y_{\RB,L}^k;\mu) =\frac{1}{\Delta t_k} \MspaceL \\
    &\hspace{3mm}+\frac{1}{2} \big(c(\mu) \Vspace+a(\mu)\mathsf\Psi^\top_\ell \Mlump \Phimat (\Pmat^\top \Phimat)^{-1} \Theta\big(\Pmat^\top \RBmat \yfec_{\RB,L}^k\big)\Pmat^\top \RBmat\big) \in \mathbb R^{\ell \times \ell}
\end{align*}
with $\Theta(\mathsf v) = \max\{0,\mathsf v\}\in\mathbb R^L$ for $\mathsf v \in \mathbb R^L$. Problem \eqref{RootRB} is now independent of the spatial FE dimension $N$. Note that we cannot guarantee the uniqueness of $\yfec_{\RB, L}^k$ and the unique solvability of the Newton iteration by the same arguments as for the FE and RB problem. In practice, a sufficiently accurate DEIM approximation of the nonlinearity will also restore the monotonicity of the max term and thus the monotonicity argument from Remark \ref{rem: NewtonSystem} can be applied again, cf. \cite{Hin20} for an analogous argumentation in the context of EIM.  We will only briefly comment on how $\Pmat$ and $\Phimat$ are generated:
\begin{enumerate}
    \item [1)] A discrete DEIM training set $\trainDEIM \subset \Pmad$ is chosen;
    \item [2)] Snapshots are generated as $\max\{0, \yfec_\delta^k(\mu)\}$ for $\mu \in \trainDEIM$ and $k = 1,\ldots,L$;
    \item [3)] The DEIM algorithm generates $\Pmat$ and $\Phimat$ from the snapshots.
\end{enumerate}

\subsection{RB-DEIM error estimation}

Often the additional error to the solution due to DEIM is not included in the error estimator. The problem is that DEIM approximates the nonlinearity on a discrete level. Therefore working with the DEIM approximation in the variational formulation is problematic. It is also problematic to efficiently incorporate the DEIM error later, since this would mean that also a Riesz representative must be calculated for the additional error quantity, which, to the best of our knowledge, is not possible in a parameter separable fashion. Thus we want to present an idea on how to already (partially) discretize the residual and incorporate the DEIM error early.

For an element $\varphi\in \test_\delta \cup \trial_\delta$ in the FE solution or test space, we denote the nodal values at time $t \in I$ by $\hat{\varphi}(t) \in \mathbb R^N$. Now we introduce the RB-DEIM residual:
\begin{align}
    \label{eq:RB_DEIM_residual}\calR_{\RB, L}(\phi;\mu)= \calF(\phi;\mu) - \calA_L(y_{\RB, L}(\mu), \phi;\mu) \quad \text{for all }\phi\in \test_\delta
\end{align}
with
\begin{align*}
    \calA_L(y_{\RB, L}(\mu), \phi; \mu) &= \calB(y_{\RB, L}(\mu),\phi; \mu) + \calN_L(y_{\RB, L}(\mu),\phi; \mu), \\
    \calN_L(y_{\RB, L}(\mu), \phi; \mu) &=a(\mu) \int_I {\langle \Phimat (\Pmat^\top \Phimat)^{-1} \max\{0, \Pmat^\top  \hat{y}_{\RB, L}(\mu) \},\hat{\phi} \rangle}_{\Mlump}\,\mathrm d t
\end{align*}
and $\langle\cdot\,,\cdot \rangle_{\Mlump}=\langle \Mlump \cdot\,,\cdot\rangle_2$, where $\langle\cdot\,,\cdot \rangle_2$ denotes the Euclidean norm. Moreover, $\calN_L$ represents the approximation obtained by mass lumping and DEIM.

\begin{remark}
The sequence of root finding problems $G_{\RB, L}^k$ can equivalently be derived by working with $\test_\RB$ as test space in \eqref{eq:RB_DEIM_residual} and applying a trapezoidal quadrature rule analogously to Section \ref{sec:FE_space_time}. This further justifies the idea of already incorporating DEIM in the variational formulation.\hfill$\Diamond$
\end{remark}

The idea is now to proceed analogously to Theorem \ref{thm: error_estimator_fix_new}. Before we state the main result of this section, let us introduce the symbols $\hat{\Hm} = L^2(I;\Mspace)$ and $\tilde{\Hm} = L^2(I;\Mlump)$, where $L^2(I;\Mspace)$ and  $L^2(I;\Mlump)$ denote the spaces of all (measurable) functions $\varphi \colon I\to\mathbb R^N$ satisfying 
\begin{align*}
& \|\varphi\|_{\hat{\Hm}}=\bigg(\int_I\|\varphi(t)\|_{\Mspace}^2\,\mathrm dt\bigg)^{1/2}<\infty,\\
& \|\varphi\|_{\tilde{\Hm}}=\bigg(\int_I\|\varphi(t)\|_{\Mlump}^2\,\mathrm dt\bigg)^{1/2}<\infty,
\end{align*}
respectively. Note that
\begin{align*}
    {\| \varphi \|}_\Hm^2={\| \hat{\varphi} \|}_{\hat{\Hm}}^2 \leq \ \| \hat{\varphi} \|^2_{\tilde{\Hm}} + \epsilon_{\tilde{\Hm}}(\varphi, \varphi),
\end{align*}
with the mass lumping error defined as
\begin{align}
    \label{MassLumpErr}
    \epsilon_{\tilde{\Hm}}(\varphi, \psi) = \vert\langle \varphi, \psi \rangle_\Hm - {\langle \hat{\varphi}, \hat{\psi} \rangle}_{\hat{\Hm}}\vert\quad \text{for all } \varphi, \psi \in \test_\delta \cup \trial_\delta.
\end{align}

\begin{proposition}
    \label{prop: error_estimator_RB_DEIM}
    Let Assumption~\emph{\ref{Assumption1}} hold. Suppose that the projection operator $\calP^\delta$ satisfies Assumption~\emph{ \ref{AssumptionProj}}. For $\mu \in \Pmad$ denote the RB-DEIM error as $e_{\RB, L}(\mu) = y_\delta(\mu) - y_{\RB, L}(\mu)$ satisfying $\calP^\delta e_{\RB, L}(\mu) \neq 0$. Then,
    \begin{align*}
        {\|e_{\RB, L}(\mu) \|}_{\test} \leq \Delta_{\RB, L}(\mu) + \Delta_{\calP^\delta}(\mu) + \Delta_L(\mu) + \Delta_{\tilde{\Hm}}(\mu),
    \end{align*}
    with DEIM approximated RB error estimator $\Delta_{\RB, L}$, projection error estimator $\Delta_{\calP^\delta}$, DEIM error estimator $\Delta_L$ and mass lumping error estimator $\Delta_{\tilde{\Hm}}$ given by
    \begin{align*}
        \Delta_{\RB, L}(\mu) & = \frac{1}{c(\mu)}\,{\|\calR_{\RB, L}(\cdot\,; \mu)\|}_{\Ym_\delta'},\\
        \Delta_{\calP^\delta}(\mu) & =\frac{a(\mu)}{c(\mu)}\,\frac{{\|\max\{0, y_\delta(\mu)\} - \max\{0, y_{\RB, L}(\mu)\}\|}_\Hm}{\|\calP^\delta e_{\RB, L}(\mu)\|_\Ym} {\|e_{\RB, L}(\mu) - \calP^\delta e_{\RB, L}(\mu)\|}_\Hm,\\
        \Delta_L(\mu) &=\frac{C_P a(\mu)}{c(\mu)} {\|\Phimat (\Pmat^\top \Phimat)^{-1}\max\{0, \Pmat^\top \hat{y}_{\RB, L}(\mu)\} - \max\{0, \hat{y}_{\RB ,L}(\mu)\}\|}_{\tilde{\Hm}},\\
        \Delta_{\tilde{\Hm}}(\mu) &= \frac{a(\mu)}{c(\mu){\|\calP^\delta e_{\RB, L}(\mu)\|}_\Ym} \Big[ \epsilon_{\tilde{\Hm}}\left(\max\left\{0, y_{\RB, L}(\mu)\right\}, \calP^\delta e_{\RB, L}(\mu)\right) \\
        &\quad + {\|\Phimat (\Pmat^\top \Phimat)^{-1}\max\{0, \Pmat^\top \hat{y}_{\RB, L}(\mu)\} - \max\{0, \hat{y}_{\RB ,L}(\mu)\}\|}_{\tilde{\Hm}} \\
        & \qquad \cdot \epsilon_{\tilde{\Hm}}(\calP^\delta e_{\RB, L}(\mu), \calP^\delta e_{\RB, L}(\mu))\Big].
    \end{align*}
\end{proposition}

\begin{proof}
    Let $\mu \in \Pmad$ be fixed. Analogously to the proof of Theorem~\ref{thm: error_estimator_fix_new} we obtain
    \begin{align*}
        & c(\mu)\,{\|\calP^\delta e_{\RB, L}(\mu)\|}_\test^2\\
        &\leq \calR_{\RB, L}(\calP^\delta e_{\RB, L}(\mu);\mu)\\
        &\quad+ a(\mu) \Big[- {\langle \max\{0, y_{\RB, L}(\mu), \calP^\delta e_{\RB, L}(\mu)\rangle}_\Hm\\
        & \hspace{15mm}+ {\langle \max\{0, y_\delta(\mu)\} -  \max\{0, y_{\RB, L}(\mu)\}, e_{\RB, L}(\mu) - \calP^\delta e_{\RB, L}(\mu)\rangle}_\Hm\\
        & \hspace{15mm}+ \langle \Phimat (\Pmat^\top \Phimat)^{-1} \max\{0, \Pmat^\top \hat{y}_{\RB, L}(\mu)\}, \calP^\delta \hat{e}_{\RB, L}(\mu) \rangle_{\tilde{\Hm}}\Big]\\
        & \leq {\|\calR_{\RB, L}(\cdot\,; \mu)\|}_{\test_\delta^\prime} {\|\calP^\delta e_{\RB, L}(\mu)\|}_\test\\
        &\quad+ a(\mu) {\|\max\{0, y_\delta(\mu)\} - \max\{0, y_{\RB, L}(\mu)\}\|}_\Hm{\|e_{\RB, L}(\mu) - \calP^\delta e_{\RB, L}(\mu)\|}_\Hm\\
        &\quad+ a(\mu) \Big\vert{\langle \Phimat (\Pmat^\top \Phimat)^{-1} \max\{0, \Pmat^\top \hat{y}_{\RB, L}(\mu)\}, \calP^\delta \hat{e}_{\RB, L}(\mu) \rangle}_{\tilde{\Hm}}\\
        & \hspace{15mm}-{\langle \max\{0,y_{\RB, L}(\mu), \calP^\delta e_{\RB, L}(\mu)\rangle}_\Hm\Big\vert.
    \end{align*}
    Now we use the mass lumping error defined in \eqref{MassLumpErr} to estimate
    \begin{align*}
        & a(\mu) \Big\vert{\langle \Phimat (\Pmat^\top \Phimat)^{-1} \max\{0, \Pmat^\top \hat{y}_{\RB, L}(\mu)\}, \calP^\delta \hat{e}_{\RB, L}(\mu) \rangle}_{\tilde{\Hm}}\\
        &\hspace{7mm}-{\langle\max\{0, y_{\RB, L}(\mu), \calP^\delta e_{\RB, L}(\mu)\rangle}_\Hm\Big\vert\\
        & \leq a(\mu) \Big\vert{\langle \Phimat (\Pmat^\top \Phimat)^{-1}\max\{0, \Pmat^\top \hat{y}_{\RB, L}(\mu)\}, \calP^\delta \hat{e}_{\RB, L}(\mu)\rangle}_{\tilde{\Hm}}\\
        &\hspace{12mm}-{\langle \max\{0, \hat{y}_{\RB, L}(\mu)\}, \calP^\delta \hat{e}_{\RB, L}(\mu)\rangle}_{\tilde{\Hm}} + \epsilon_{\tilde{\Hm}}(\max\{0, y_{\RB, L}(\mu)\}, \calP^\delta e_{\ell, L}(\mu)\Big\vert\\
        & \leq a(\mu)\Big[{\|\Phimat (\Pmat^\top \Phimat)^{-1}\max\{0, \Pmat^\top \hat{y}_{\RB, L}(\mu)\} - \max\{0, \hat{y}_{\RB ,L}(\mu)\}\|}_{\tilde{\Hm}}{\|\calP^\delta \hat{e}_{\RB, L}(\mu)\|}_{\tilde{\Hm}}\\
        &\hspace{12mm}+ \epsilon_{\tilde{\Hm}}\big(\max\{0, y_{\RB, L}(\mu)\}, \calP^\delta e_{\RB, L}(\mu)\big)\Big]\\
        & \leq a(\mu) \Big[{\|\Phimat (\Pmat^\top \Phimat)^{-1}\max\{0, \Pmat^\top \hat{y}_{\RB, L}(\mu)\} - \max\{0, \hat{y}_{\RB ,L}(\mu)\}\|}_{\tilde{\Hm}}\\
        &\hspace{12mm}\cdot\big[ c_P\,
        {\|\calP^\delta e_{\RB, L}(\mu))\|}_{\test}+ \epsilon_{\tilde{\Hm}}\big(\calP^\delta e_{\RB, L}(\mu), \calP^\delta e_{\RB, L}(\mu)\big)\big]\\
        &\hspace{12mm}+ \epsilon_{\tilde{\Hm}}\big(\max\{0, y_{\RB, L}(\mu)\}, \calP^\delta e_{\RB, L}(\mu)\big)\Big],
    \end{align*}
    where we have used
    \begin{align*}
        \|\calP^\delta \hat{e}_{\RB, L}(\mu)\|_{\tilde{\Hm}} \leq \|\calP^\delta \hat{e}_{\RB, L}(\mu)\|_{\Hm} + \epsilon_{\tilde{\Hm}}\big(\calP^\delta e_{\RB, L}(\mu), \calP^\delta e_{\RB, L}(\mu)\big)
    \end{align*}
    and Poincaré's inequality with constant $c_P > 0$ in the last estimate. All in all we obtain
    \begin{align*}
        &{\|\calP^\delta e_{\RB, L}(\mu)\|}_\test \leq \frac{1}{c(\mu)}\,{\|\calR_{\RB, L}(\cdot\,; \mu)\|}_{\test_\delta'} \\
        & + \frac{a(\mu)}{c(\mu)} \frac{{\|\max\{0, y_\delta(\mu)\} - \max\{0, y_{\RB, L}(\mu)\}\|}_\Hm}{\|\calP^\delta e_{\RB, L}(\mu)\|_\test}\,{\|e_{\RB, L}(\mu) - \calP^\delta e_{\RB, L}(\mu)\|}_\Hm\\
        & + \frac{a(\mu)c_P}{c(\mu)} \|\Phimat (\Pmat^\top \Phimat)^{-1}\max\{0, \Pmat^\top \hat{y}_{\RB, L}(\mu)\} - \max\{0, \hat{y}_{\RB ,L}(\mu)\}\|_{\tilde{\Hm}}\\ 
        &  + \frac{a(\mu)}{c(\mu)\|\calP^\delta e_{\RB, L}(\mu)\|_\test} \Big[ \|\Phimat (\Pmat^\top \Phimat)^{-1}\max\{0, \Pmat^\top \hat{y}_{\ell, L}(\mu)\} - \max\{0, \hat{y}_{\ell ,L}(\mu)\}\|_{\tilde{\Hm}}\\
        & \qquad \cdot \epsilon_{\tilde{\Hm}}(\calP^\delta e_{\ell, L}(\mu), \calP^\delta e_{\ell, L}(\mu)) + \epsilon_{\tilde{\Hm}}\big(\max\left\{0, y_{\RB, L}(\mu)\right\}, \calP^\delta e_{\RB, L}(\mu)\big)\Big]
    \end{align*}
    All in all the desired estimate follows.
\end{proof}

\begin{remark}
    \label{rem: numerical_projection_error}
    \begin{enumerate}
        \item [a)] Note that analogously to $\Delta_{\calP^\delta}$ also $\Delta_{\tilde{\Hm}}$ can be made small independently of the chosen RB space, by choosing a sufficiently rich spatial FE space. This follows as as a consequence of \cite[Lemma~1]{Z12} adopted to the present PDE. Therefore we will use $\Delta_{\RB, L} + \Delta_L$ as error estimator in our numerical experiments.
        \item [b)] The computation of $\Delta_L$ requires an evaluation of the nonsmoothness on FE level. This cannot be avoided, but is computationally still cheap compared to the other computations.\hfill$\Diamond$
\end{enumerate}
\end{remark}

\subsection{Adaptive DEIM}

In Proposition \ref{prop: error_estimator_RB_DEIM} a decomposition of the total error into a temporal mesh-dependent part $\Delta_{\mathcal{P}^\delta}$, a spatial mesh-dependent part $\Delta_{\tilde{\Hm}}$, an RB dependent part $\Delta_{\RB, L}$ and a DEIM dependent part $\Delta_L$ has been proposed. As already argued, we will neglect the mesh-dependent parts, since they can be chosen arbitrarily small by sufficiently rich spatial and temporal FE spaces. Now typically a DEIM basis would be generated in an offline phase to efficiently approximate the nonsmoothness in a lower dimensional subspace as described in Section \ref{sec: ClassicalDEIM}. This process is decoupled from the RB basis generation described in Section \ref{sec: RB_basis}. More recently ideas have been presented on how to combine RB and DEIM basis generation into one adaptive process, see e.g. \cite{Che20,Pet15}. Advantages of such ideas are control over the DEIM error in the offline phase and thus a DEIM basis of suitable size. Also the combination of RB and DEIM basis generation leads to fewer computational costs. The adaptive DEIM algorithm in \cite{Pet15} is taylored towards optimization and updates the DEIM basis without increasing its size. This results in only locally good approximation quality, which in optimization is sufficient. In \cite{Che20} the authors essentially calculate a completely new DEIM basis whenever it needs to be updated. On one hand, this is cheap, since no additional snapshot generation is necessary and the size of the singular value decomposition necessary for DEIM is usually still small compared to classical offline DEIM. On the other hand, this usually means that all preassembled DEIM data must be thrown away whenever the DEIM basis is updated. Therefore we suggest a successive DEIM approach for the space-time setting, which might result in a (slightly) larger DEIM basis, but allows to reuse previous information. We summarize this approach in Algorithm \ref{algo: adaptive_POD_Greedy_RB_DEIM}.

\begin{algorithm}[h!]
    \caption{(Adaptive POD-Greedy RB-DEIM Method)}\label{algo: adaptive_POD_Greedy_RB_DEIM}
    \begin{algorithmic}[1]
        \Require{Initial DEIM data $\Pmat_0$, $\Phimat_0$, $\Imat_0$, discrete training set of parameters $\train \subset\Pmad$, error tolerances $\RBtol, \DEIMtol > 0$ and $\Totaltol \geq \RBtol + \DEIMtol$;}
        \State{Set $\ell = 0$, $\mathsf{\Psi}_{0} = \emptyset$, $\Vs_{0} =  \{0\}$, $\Pmat = \Pmat_0$, $\Phimat = \Phimat_0$ and $\Imat= \Imat_0$;}
		\State{Compute $\tilde{\mu} \in \argmax\{\Delta_{\RB, L}(\mu) + \Delta_L(\mu)\,\vert\,\mu \in \train\}$;}
		\State{Set $\varepsilon_1=\Delta_{\RB, L}(\tilde{\mu})$ and $\varepsilon_2=\Delta_L(\tilde{\mu})$;}
        \While{$\varepsilon_1  + \varepsilon_2  > \Totaltol$}
            \State{\em \% RB part;}
            \If{$\varepsilon_1 > \RBtol$}
                \State{Set $e_{\ell + 1}^k = y_\delta^k(\tilde{\mu}) - \calP_{\Vs_\ell} y_\delta^k(\tilde{\mu})$ for $k = 1,\ldots,K$;}
                \State{Compute $\psi_{\ell + 1} \in \POD_1(\{e_{\ell + 1}^k\}_{k = 1}^K)$;}
			    \State{Set $\mathsf{\Psi}_{\ell + 1} = \RBmat \cup \{\psi_{\ell + 1}\}$, $V_{\ell + 1} = V_\ell \oplus \text{span}(\psi_{\ell + 1})$ and $\ell = \ell + 1$;}
        \EndIf   
		\State{\em \% DEIM part;}
		\State{Set $L=\lfloor \log_{10}(\varepsilon_2/\DEIMtol) \rfloor$;}
		\If{$L = 0$}
            \State{$L=1$ (avoid stagnation);}
        \EndIf
        \State{Set $E = [\max\{0, y_\delta^1(\tilde{\mu}\}, \ldots, \max\{0, y_\delta^K(\tilde{\mu}\}]$;\hfill{\em \% snapshot matrix}}
		\State{Put $E=E-\Phimat\Phimat^\top E$;\hfill{\em \% subtract previous DEIM information}}
        \State{Set $\tilde L=\min\{L,\mathrm{rank}(E)\}$;}
		\State{Compute POD basis $\tilde{\Phimat}$ of length $\tilde{L}$ w.r.t. identity matrix for $E$;}
		\If{$\Phimat = []$ (no initial DEIM basis)}
            \State{Compute $i \in \argmax_{k = 1,\ldots,N} \vert(\tilde{\phi}_1)_k\vert$;}
			\State{Set $\Phimat = [\tilde{\phi}_1]$, $\Imat = [i]$, $\Pmat = [e_i]$ and $\tilde{j} = 2$;}
		\Else
			\State{Set $\tilde{j} = 1$;}
        \EndIf
		\For{$j = \tilde{j}, \ldots, \tilde{L}$} 			
			\State{Solve $(\Pmat^\top \Phimat) \gamma = \Pmat^\top \tilde{\phi}_j$;}
			\State{Set $r = \tilde{\phi}_j - \Phimat \gamma$;}
			\State{Compute $i\in\argmax_{k=1,\ldots,N}\vert r_k\vert$;}
			\State{Set $\Phimat = [\Phimat, \tilde{\phi}_j]$, $\Imat = [\Imat, i]$ and $\Pmat = [\Pmat, e_i]$;}
		\EndFor	
		\State{\em\% Error update;}
		\State{Compute $\tilde{\mu} \in \argmax\{\Delta_{\RB, L}(\mu) + \Delta_L(\mu)\,\vert\,\mu \in \train\}$;}
		\State{Set $\varepsilon_1=\Delta_{\RB, L}(\tilde{\mu})$ and $\varepsilon_L=\Delta_L(\tilde{\mu})$;}
        \EndWhile\hfill\\
        \Return{Reduced basis $\RBmat$, spatial RB space $\Vs_\ell$, DEIM data $\Pmat$, $\Phimat$, $\Imat$.}
    \end{algorithmic}
\end{algorithm}

To avoid stagnation at the end of the algorithm one should e.g. choose $\RBtol = 0.1 \Totaltol$ and $\DEIMtol = 0.01 \Totaltol$. This ensures that the algorithm aims at individual tolerances $\RBtol$ and $\DEIMtol$ well below $\Totaltol$. Also it should be ensured that the DEIM basis is always at least as large as the RB basis for numerical stability. This is the reason why there is the possibility to start the algorithm with an initial DEIM basis. For more details we refer to \cite{Che20}. The algorithm performs best, if the construction of a RB and a DEIM basis are of equal difficulty. This is often not the case and usually the construction of a DEIM basis is more difficult, as also our numerical experiments in Section \ref{sec: numerics} suggest. Heuristics to let the DEIM basis grow faster exist, e.g. $\log_{10}$ could be replaced by $\log_{2}$ in line 11. This can be very problem specific and is beyond the scope of this work.
%
%
In our implementation a maximum of one element is added to the RB basis per iteration. In \cite{Che20} a logarithmic term analogue to the one for the DEIM basis is used. In our opinion this might lead to too many elements being added during the first iterations and wrong elements removed later on. But of course the algorithm could also be extended in this fashion. 

\section{Numerical experiments}
\label{sec: numerics}

In this section we present two numerical examples to investigate the performance of the different space-time RB approaches. Parameter choices fixed for both examples are summarized in Table \ref{tab: parameters}.
\begin{table}[!ht]
    \begin{center}
    \begin{tabular}{c|cc|ccc}\toprule
        PDE & \multicolumn{2}{c}{discretization} & \multicolumn{3}{|c}{tolerances} \\\midrule
        $\Omega$ & $K$ & $1/h$ & $\Totaltol$ & $\RBtol$& $\DEIMtol$\\ \midrule
        $(0,1)^2$ & $400$ & $\{50, 100, 200\}$ & $10^{-3}$ & $10^{-4}$ & $10^{-5}$ \\\bottomrule
    \end{tabular}\\[2mm]
    \caption{Fixed parameters for the numerical examples. Recall that $K$ and $h$ denote the number of time grid points and the spatial step-size, respectively.} 
    \label{tab: parameters}
    \end{center}
\end{table}
\newline
The adaptive DEIM algorithm uses an initial DEIM basis of size two generated from one snapshot.

Our code is implemented in Python 3 and uses FEniCS (see \cite{FENICS}) for the matrix assembly. Sparse memory management and computations are implemented with SciPy (see \cite{Scipy}). All computations below were run on an Ubuntu 20.04 notebook with 32 GB main memory and an Intel Core i7-8565U CPU.

\subsection{Example 1}

We choose $T = 20$ for the maximal time horizon and set $\mathscr P_\mathsf{ad} = [-10,10]$, $a(\mu) = 1 + 2\,\vert\mu\vert$ and $c(\mu) = 5 / (5+\vert\mu\vert)$. Furthermore the right-hand sides are
%
\begin{align*}
    f(t;\mu)(x_1, x_2) = \underbrace{10 \sin\Big(\frac{4\pi t}{T}\Big) \sqrt{1 + t}}_{\text{time}}\underbrace{\Big(\frac{1}{2} - x_1\Big) \sin(\pi x_1) \sin(\pi x_2)}_{\text{space}} \underbrace{\mu.}_{\text{parameter}}
\end{align*}
For the classical RB and RB-DEIM approach we fix a training set of $60$ equidistant parameters and a test set of $100$ equidistant parameters. For DEIM we additionally choose the same training set as for the RB offline phase and choose the same parameter $L$ we obtain from the adaptive RB-DEIM. Note that this means that already information from the adaptive approach is used in the classical RB-DEIM approach, but it allows for a better comparison of the results, which are shown in Table \ref{tab: ex1}. 
\begin{table}[!h]
    \begin{center}
        {\tiny\begin{tabular}{ccccccccc}\toprule
        $1/h$ & av. & time & av. & av. & av. proj. & av. & size& size\\
         & sp.-up & offline & error & est. & error & eff. & RB & DEIM\\\midrule
        \multicolumn{8}{c}{RB (true error)}\\ \midrule
        \phantom{1}50  & \phantom{14}5.95 & $7.76 \cdot 10^2$ & $2.73 \cdot 10^{-4}$ &  $2.88 \cdot 10^{-4}$ & $4.73 \cdot 10^{-6}$ & 1.04 & 15 & -- \\
        100 & \phantom{1}39.54 & $4.54 \cdot 10^3$ & $5.14 \cdot 10^{-4}$ & $5.17 \cdot 10^{-4}$ & $1.72 \cdot 10^{-6}$ & 1.00 & 16 & -- \\
        200 & \phantom{1}93.72 & $2.03 \cdot 10^4$ & $5.62 \cdot 10^{-4}$ & $5.64 \cdot 10^{-4}$ & $1.44 \cdot 10^{-6}$ & 1.00 & 15 & -- \\\midrule
        \multicolumn{8}{c}{RB (estimator)}\\ \midrule
        \phantom{1}50  & \phantom{14}5.48 & $9.39 \cdot 10^2$ & $2.72 \cdot 10^{-4}$ & $2.88 \cdot 10^{-4}$ & $4.74 \cdot 10^{-6}$ & 1.04 & 15 & --  \\
        100 & \phantom{1}39.66 & $3.54 \cdot 10^3$ & $6.62 \cdot 10^{-4}$ & $6.65 \cdot 10^{-4}$ & $1.51 \cdot 10^{-6}$ & 1.00 & 15 & --\\
        200 & \phantom{1}91.17 & $1.32 \cdot 10^4$ & $6.21 \cdot 10^{-4}$ & $6.23 \cdot 10^{-4}$ & $1.34 \cdot 10^{-6}$ & 1.00 & 15 & --\\ \midrule
        \multicolumn{8}{c}{RB-DEIM (estimator)}\\ \midrule
        \phantom{1}50  & \phantom{14}8.06 & $3.07 \cdot 10^2 + 5.67 \cdot 10^2$ & $1.86 \cdot 10^{-4}$ & $3.22 \cdot 10^{-4}$ & $1.56 \cdot 10^{-6}$ & 1.59 & 18 & 72 \\
        100 & \phantom{1}94.29 & $3.51 \cdot 10^3 + 1.86 \cdot 10^3$ & $1.78 \cdot 10^{-4}$ & $3.80 \cdot 10^{-4}$ & $3.28 \cdot 10^{-7}$ & 3.37 & 27 & 66\\
        200 & 879.23 & $1.82 \cdot 10^4 + 6.29 \cdot 10^3$ & $4.28 \cdot 10^{-4}$ & $6.19 \cdot 10^{-4}$ & $3.54 \cdot 10^{-7}$ & 1.45 & 19 & 87 \\ \midrule
        \multicolumn{8}{c}{RB-DEIM (estimator, adaptive)}\\ \midrule
        \phantom{1}50  &  \phantom{14}7.68 & $6.43 \cdot 10^2$ & $1.66 \cdot 10^{-4}$ & $2.93 \cdot 10^{-4}$ & $8.90 \cdot 10^{-7}$ & 1.84 & 20 & 72 \\
        100 & \phantom{1}96.84 & $1.67 \cdot 10^3$ & $2.65 \cdot 10^{-4}$ & $5.10 \cdot 10^{-4}$ & $1.11 \cdot 10^{-6}$  & 2.43 & 19 & 66 \\
        200 & 882.92 & $6.62 \cdot 10^3$ & $2.72 \cdot 10^{-4}$ & $5.10 \cdot 10^{-4}$ & $9.77 \cdot 10^{-7}$ & 2.33 & 19 & 87 \\\bottomrule
        \end{tabular}\\[2mm]}
        \caption{{\sf{Example 1}.} Comparison of RB and RB-DEIM approaches for different spatial discretizations $1/h$. For RB-DEIM (estimator) offline time is given as DEIM + RB time.}
        \label{tab: ex1}
    \end{center}
\end{table}

First of all, we can observe mesh independence in the average efficiency and the size of the RB basis for the approaches without DEIM. If DEIM is used, the results are slightly mesh dependent, due to the mesh dependency of DEIM. Also the average error and error estimator are in the same regions and far below the tolerance of $10^{-3}$ for all approaches and independently of the mesh. We can observe that the RB error estimator leads to slight decreases in offline computational time, but only the introduction of DEIM allows for efficient evaluations of the estimator. Unfortunately the offline computational cost of DEIM outweighs this advantage. Finally we can observe that the adaptive RB-DEIM approach leads to the largest reduction in offline computational cost, whilst simultaneously giving equivalent benefits as classical RB-DEIM in average speed-up in the online phase. Especially on finer grids the advantages of an RB-DEIM approach are evident. 
Last but not least, we want no mention two noticable entries in Table \ref{tab: ex1}. First, the RB basis basis without estimator contains one more element than that with estimator for $1/h = 100$. This is atypical, but of course possible due to the outstanding efficiency and the fact that the greedy procedure is only locally choosing the best update in every iteration. Second, the RB-DEIM basis without adaptivity is rather large and contains $27$ elements for $1/h = 100$. This comes from the fact that the procedure struggles to reach the desired overall tolerance in this case.

\begin{figure}[!ht]
    \begin{center}
        \subfigure
	       {\includegraphics[width=0.47\textwidth]{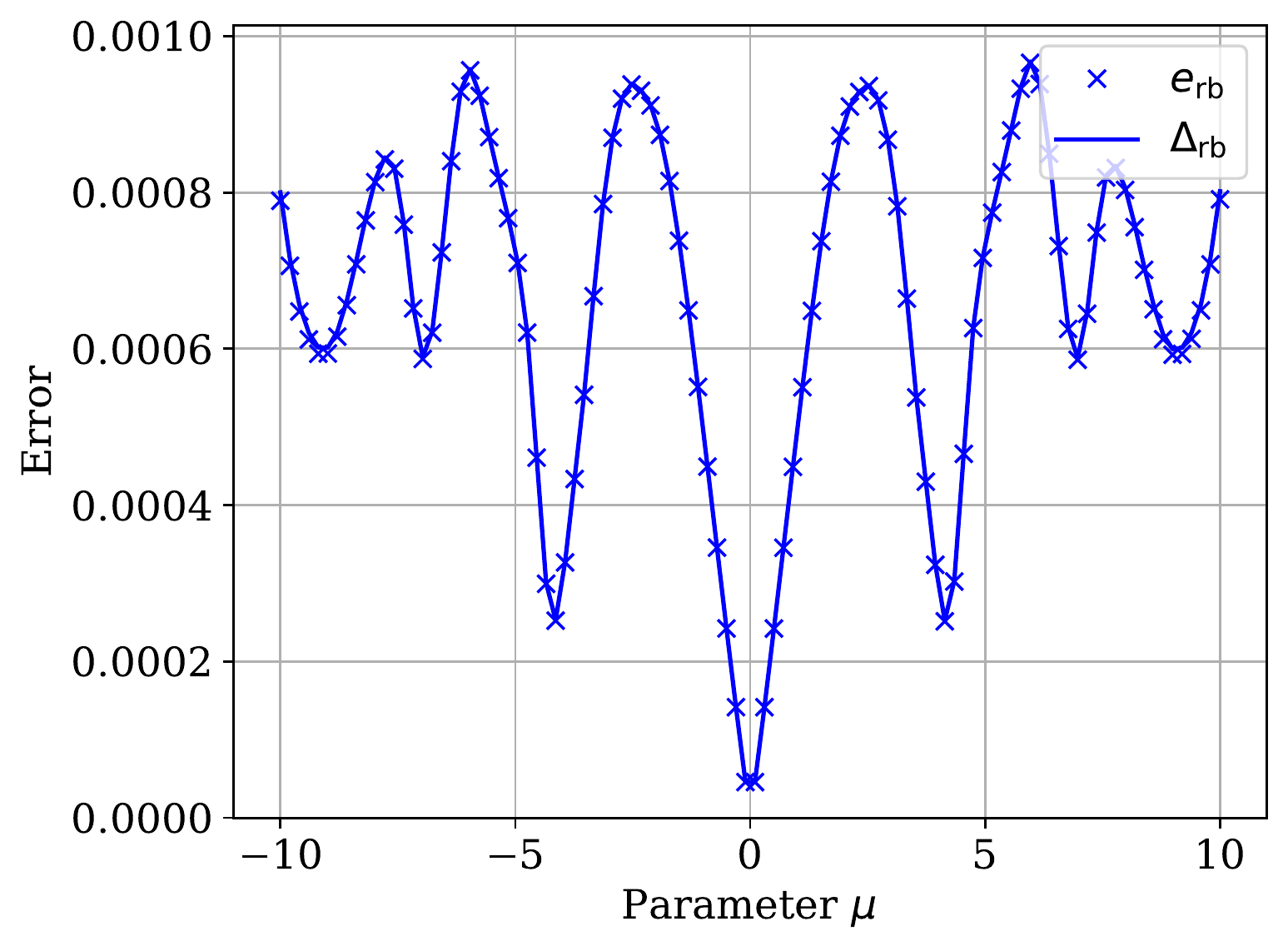}}
        \subfigure
	       {\includegraphics[width=0.47\textwidth]{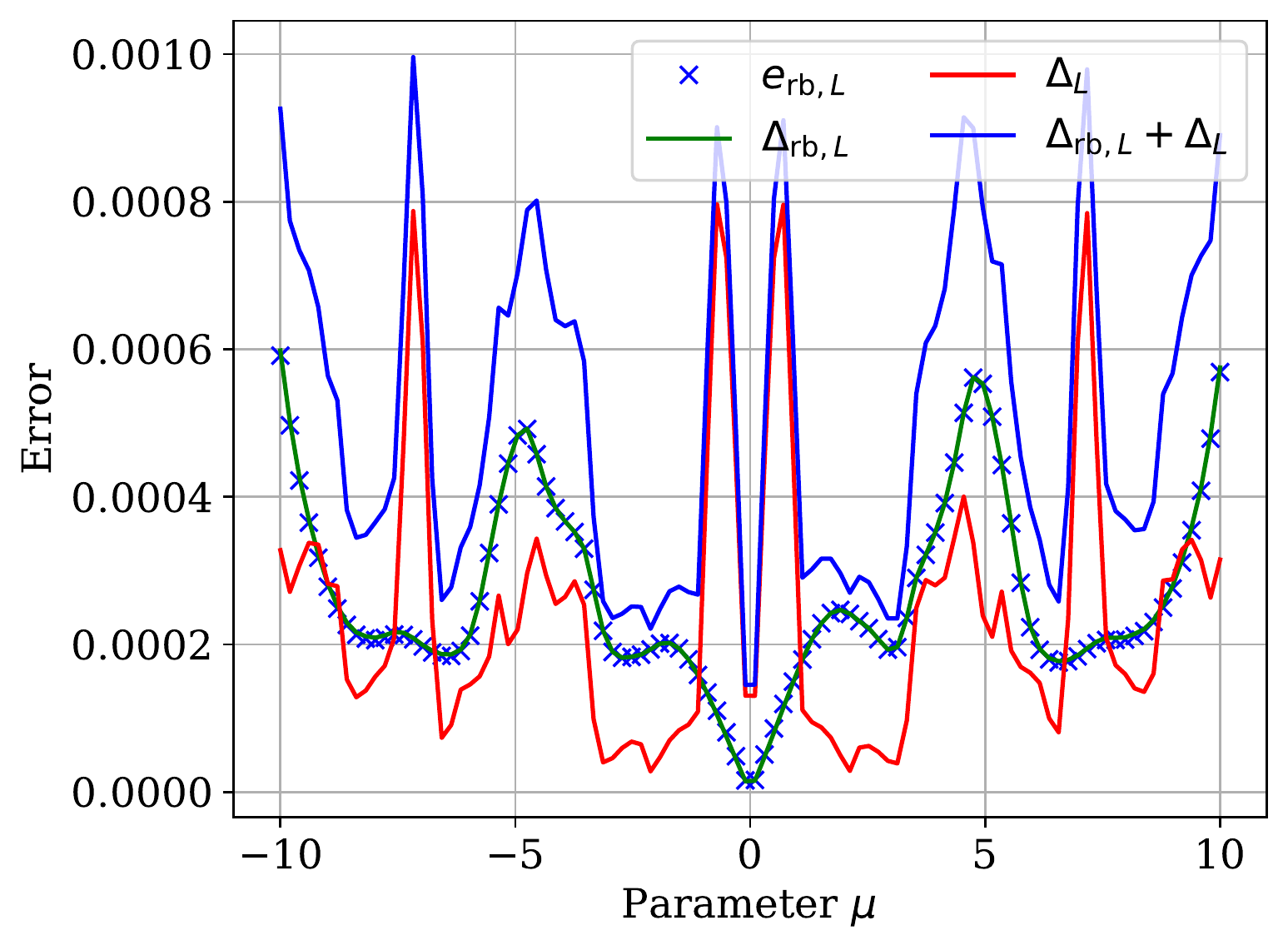}}
    \end{center}
    \caption{{\sf{Example 1}.} Error and RB error estimator (left) and error and RB-DEIM error estimator (right) with $1/h = 100$.}
    \label{fig: ex1}
\end{figure}

To further investigate the error estimator, it is shown in Figure \ref{fig: ex1} together with the error on the test set.

We can see that the RB error estimator clearly mimics the behavior of the true error. In the RB-DEIM approach, the additional DEIM error estimator is more volatile. Furthermore the additional DEIM error estimator is necessary to avoid underestimation, though this is hard to observe in the plot.

\begin{figure}[!ht]
    \begin{center}
        \subfigure
	       {\includegraphics[width=0.47\textwidth]{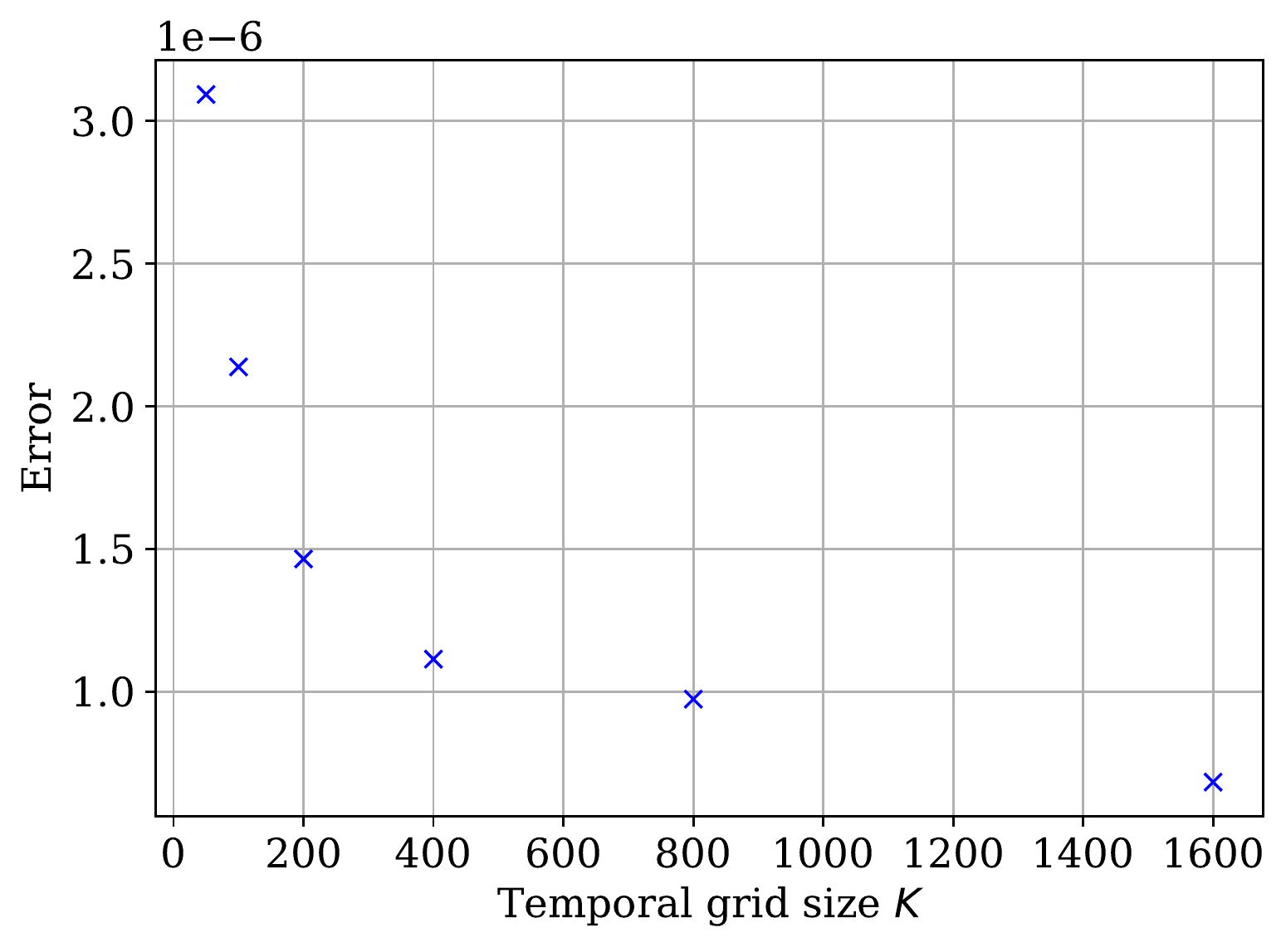}}
        \subfigure
	       {\includegraphics[width=0.47\textwidth]{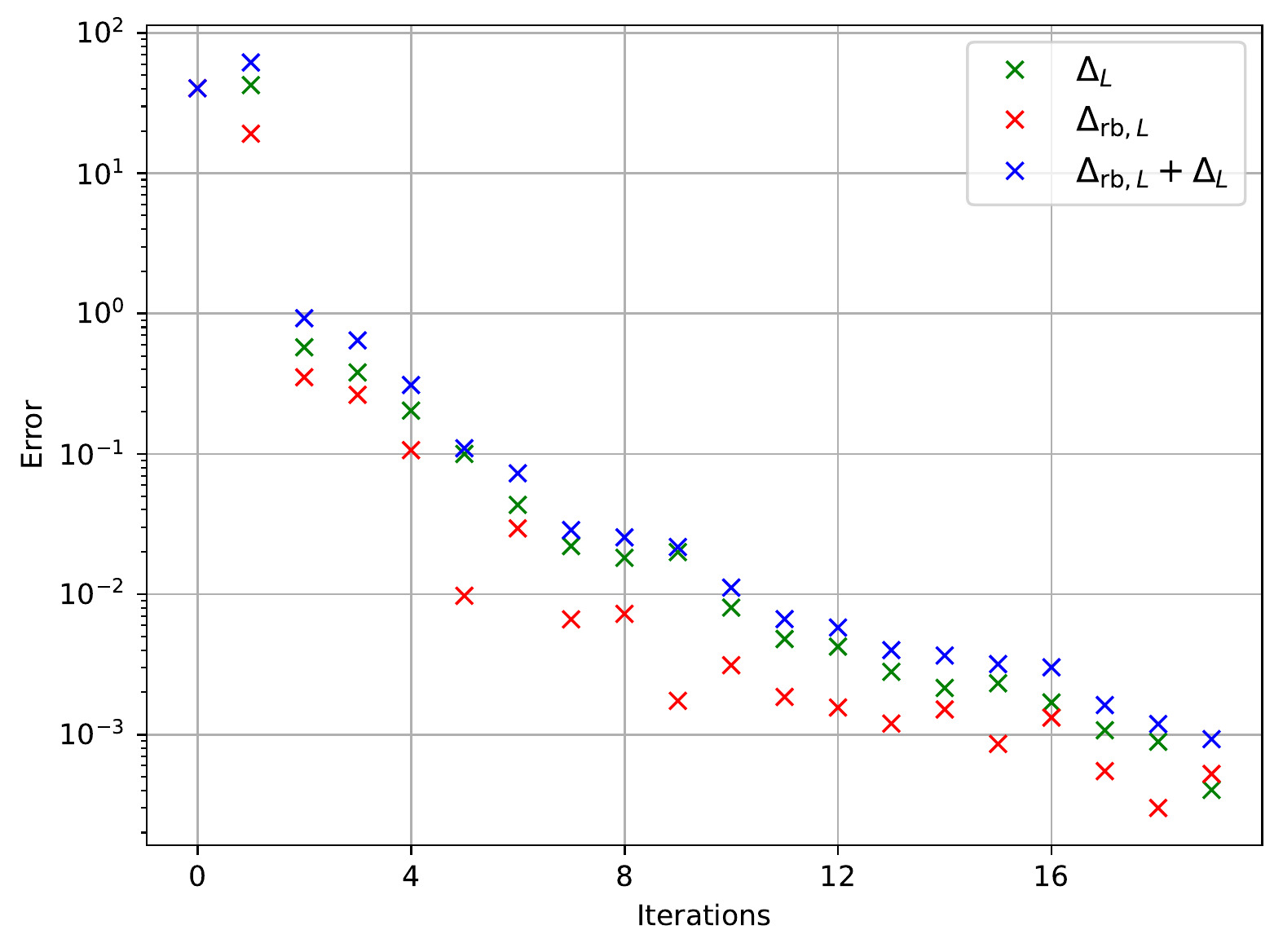}}
    \end{center}
    \caption{{\sf{Example 1}.} Left: Projection error $\Delta_{\mathcal P^\delta}$ for varying temporal grid size $K$. Right: Evolution of the total error $\Delta_{\RB, L} + \Delta_L$ and its individual components during training of the adaptive RB-DEIM. Both with $1/h = 100$.}
    \label{fig: ex1_errors}
\end{figure}

Finally we want to investigate the projection error and the evolution of the different error quantities during the adaptive RB-DEIM basis generation. The average projection error, on the previously introduced test set, is shown in Figure \ref{fig: ex1_errors} on the left and the evolution of the errors during training is shown on the right.
In the left plot we observe that the projection error $\Delta_{\mathcal P^\delta}$ decreases with increasing time grid points as proposed in Proposition~\ref{prop: time_convergence}. However, a linear convergence can only be expected asymptotically. In particular, the regime of linear convergence seems to not be reached yet. This can be expected due to the quotient in the projection error estimator, which is only asymptotically one. Nonetheless the values of the projection error estimator are very small and can be neglected as suggested. In the right plot we can see that the RB error tends to be lower than the DEIM error, but both decrease sufficiently fast.

All in all this example clearly suggests the usage of the RB error estimator over the calculation of the true error and the usage of the adaptive DEIM approach over the classical DEIM approach. We also investigated the projection error estimator and have shown that its small values justify neglection.

\subsection{Example 2}

In the second numerical example we choose $T = 10$ for the maximal time horizon and use a two-dimensional parameter space with $\mathscr P_\mathsf{ad{}} = [-2, 2]^2$, $a(\mu) = 1 + 5\,\|\mu\|_2$ and $c(\mu) = 3/(1 + \vert\mu_1\vert)$. Furthermore, the right-hand sides $f(t; \mu)$ are chosen to be
\begin{align*}
    f(t,\bx;\mu) = \underbrace{10 \sin\Big(\frac{4\pi t}{T}\Big) \sqrt{1 + t}}_{\text{time}}\,\underbrace{\left\{
    \begin{aligned}
        &x_1 x_2 \mu_1,&&\text{for } x_1 \leq \frac{1}{2},\\
        &x_1^2 x_2^2 \mu_2,&&\text{otherwise}.
    \end{aligned}
    \right.}_{\text{space and parameter}}
\end{align*}
This time we only present results for RB with error estimator and the adaptive RB-DEIM approach, since based on the theoretical results and Example 1, these are the most promising strategies. For RB we fix a training set of $144$ equidistant parameters and a test set of $225$ equidistant parameters. The results for the different spatial discretizations are shown in Table \ref{tab: ex2}.
\begin{table}[!h]
    \begin{center}
        {\tiny\begin{tabular}{ccccccccc}\toprule
            $1/h$ & av.  & time & av. & av. & av. proj. & av. & size & size\\
            & sp.-up & offline & error & est. & error & eff. & RB & DEIM \\ \midrule
            \multicolumn{8}{c}{RB (estimator)}\\ \midrule
            \phantom{1}50 & \phantom{14}5.40 & $2.08 \cdot 10^3$ & $3.09 \cdot 10^{-4}$ & $3.14 \cdot 10^{-4}$ & $5.66 \cdot 10^{-7}$ & 1.01 & 13 & -- \\
            100 & \phantom{1}45.79 & 5.98 $\cdot 10^3$ & $3.16 \cdot 10^{-4}$ & $3.22 \cdot 10^{-4}$ & $6.03 \cdot 10^{-7}$ & 1.01 & 13 & --\\
            200 & \phantom{1}96.04 & $2.82 \cdot 10^4$ & $3.15 \cdot 10^{-4}$ & $3.21 \cdot 10^{-4}$ & $6.02 \cdot 10^{-7}$ & 1.01 & 13 & -- \\\midrule
            \multicolumn{8}{c}{RB-DEIM (estimator, adaptive)}\\\midrule
            \phantom{1}50 & \phantom{14}7.58 & $3.23 \cdot 10^3$ & $7.25 \cdot 10^{-5}$ & $1.04 \cdot 10^{-3}$ & $1.96 \cdot 10^{-7}$ & 10.00 & 29 & 177 \\
            100 & 127.06 & $7.97 \cdot 10^3$ & $7.04 \cdot 10^{-5}$ & $8.42 \cdot 10^{-4}$ & $1.93 \cdot 10^{-7}$ & \phantom{1}9.63 & 30 & 212 \\
            200 & 858.93 & $2.02 \cdot 10^4$ & $5.78 \cdot 10^{-5}$ & $8.11 \cdot 10^{-4}$ & $2.25 \cdot 10^{-7}$ & 10.70 & 29 & 222\\\bottomrule
        \end{tabular}\\[2mm]}
        \caption{{\sf{Example 2}.} Comparison of RB and RB-DEIM approaches for different spatial discretizations $1/h$. For RB-DEIM (estimator) offline time is given as DEIM + RB time.}
        \label{tab: ex2}
    \end{center}
\end{table}
Again the RB results show mesh independence in the average error, error estimator and efficiency and the size of the RB basis. This time the average efficiency for RB-DEIM is significantly larger (approx. $10$ compared to approx. $1.01$) which is one of the reasons for the increased DEIM basis. Note that the RB basis is also significantly larger (approx. $30$) than for the RB approach ($13$), but both bases are mesh independent in size. This is a typical behavior that occurs since for the given example DEIM approximation is more difficult than RB approximation. Thus the RB basis is chosen in a way that a tolerance of $10^{-4}$ is reached, so that the DEIM tolerance only needs to be below $\Totaltol = 10^{-3}$. Furthermore the size of the DEIM basis is not mesh independent, but slightly increases for finer meshes.

\begin{figure}[!ht]
    \begin{center}
        \subfigure
	    {\includegraphics[width=0.47\textwidth]{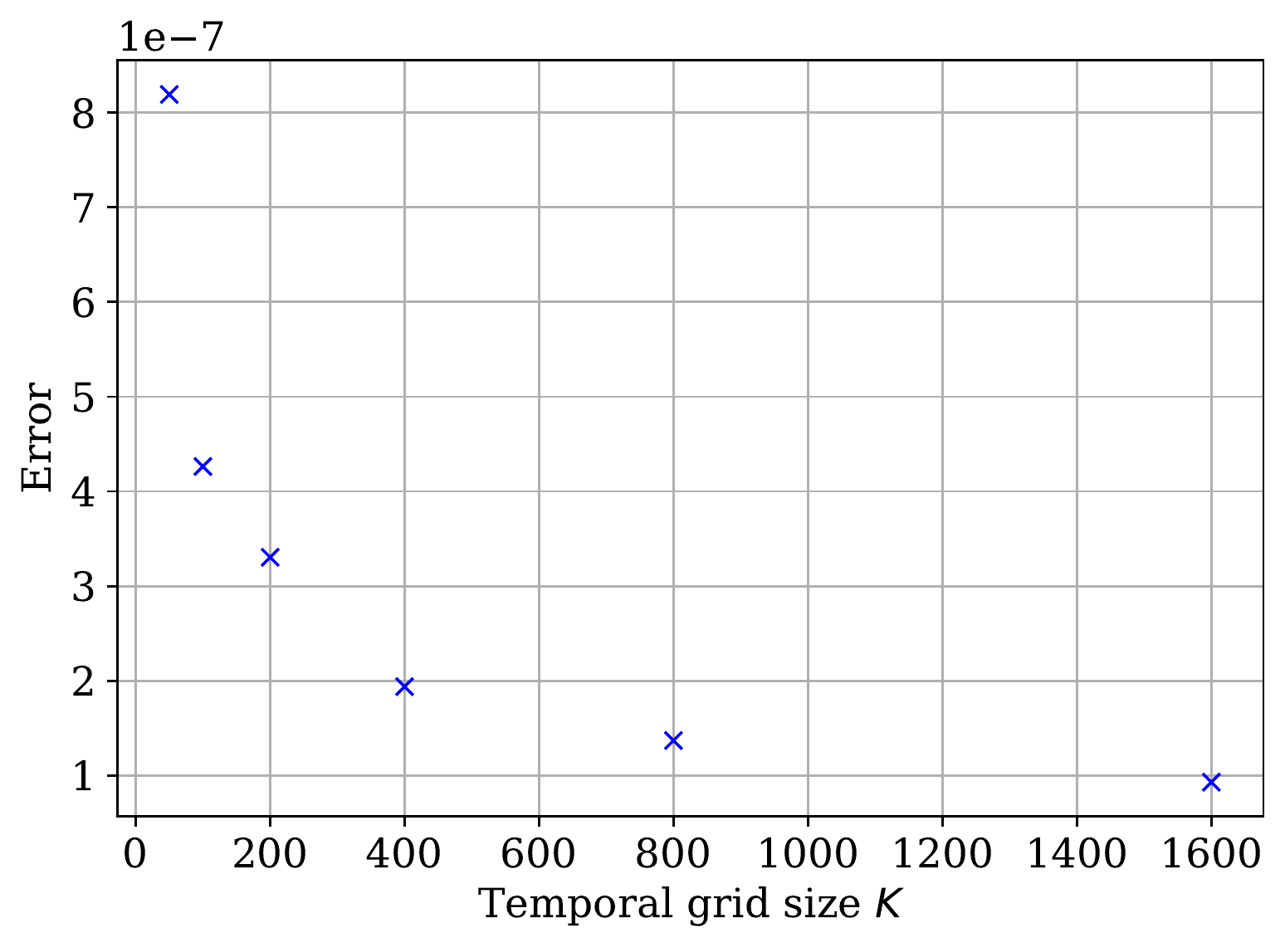}}
        \subfigure
	    {\includegraphics[width=0.47\textwidth]{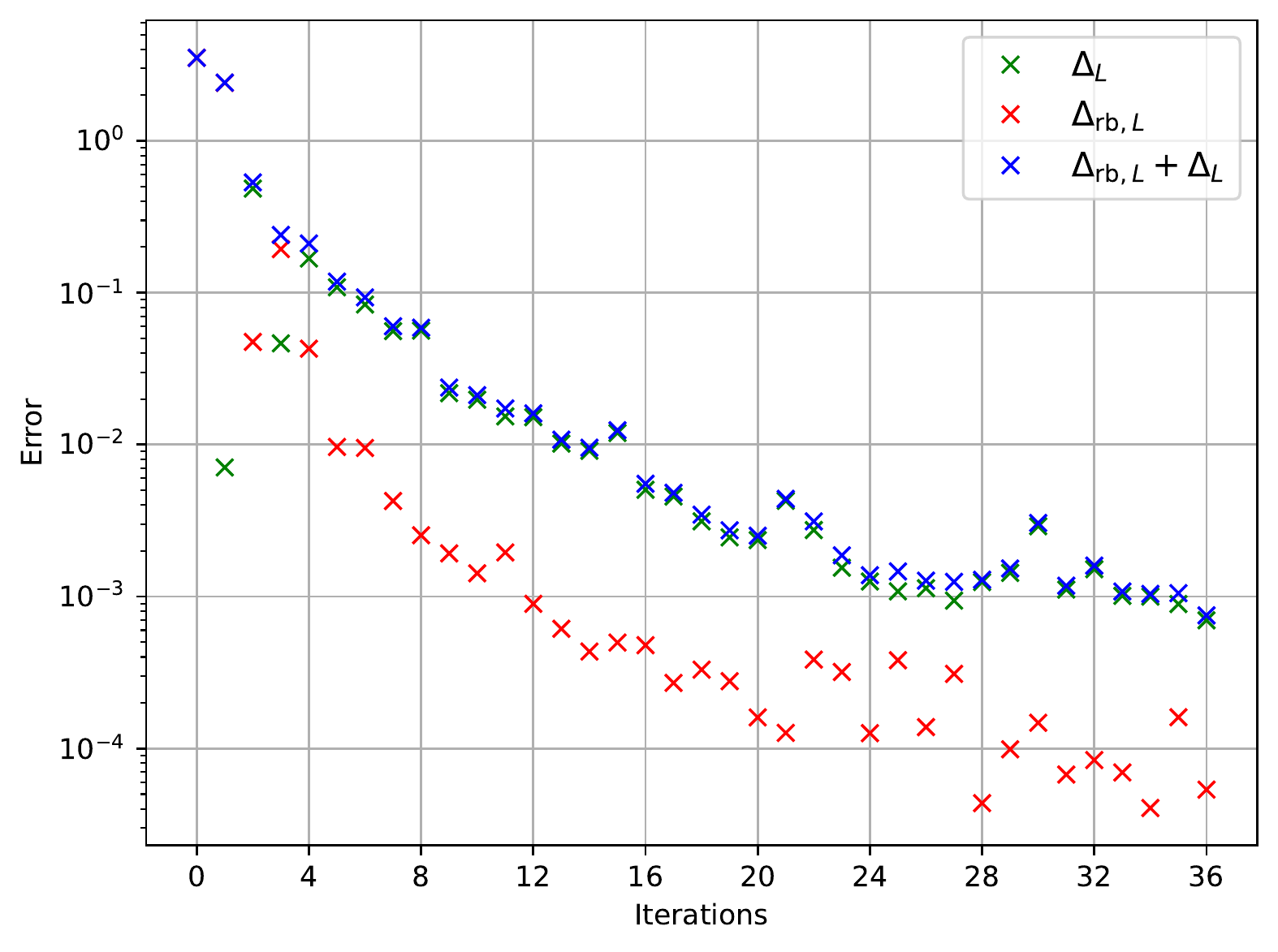}}
        \caption{{\sf{Example 2}.} Left: Projection error $\Delta_{\calP^\delta}$ for varying temporal grid size $K$. Right: Evolution of the total error $\Delta_{\RB, L} + \Delta_L$ and its individual components during training of the adaptive RB-DEIM. Both with $1/h = 100$.}
        \label{fig: ex2_errors}
    \end{center}
\end{figure}

Again we want to investigate the projection error and the evolution of the different error quantities during the adaptive RB-DEIM basis generation. The average projection error is shown in Figure \ref{fig: ex2_errors} on the left and the evolution of the errors during training is shown on the right.

In the left plot, we can observe that the average projection error is again neglectable, due to its small values and we can observe the same convergence behavior as discussed in detail for the previous example. In the right plot, we can observe that in contrast to Example 1, the slow convergence of the DEIM error slows down the overall convergence and leads to an RB error below $\RBtol = 10^{-4}$. Therefore the DEIM error in essence only needs to be below $\Totaltol = 10^{-3}$.
\newline 
The second example again suggests the usage of the RB error estimator over the calculation of the true error, due to the excellent average efficiency. Again it is unproblematic to neglect the projection error estimator. Note that we have also excluded the RB-DEIM approach without adaptivity, since no convergene can be observed with the size $L$ of the DEIM basis taken from the adaptive approach. This suggests advantages of the adaptive approach, but we can also observe that the adaptive algorithm begins to struggle when the DEIM approximation is significantly harder than the RB approximation.
\section{Conclusion}
\label{sec: conclusion}
We have introduced a novel space-time a-posteriori error estimator for a nonsmooth parabolic PDE. The numerical results show promising speed-up and good efficiency. Nonetheless the model order reduction is limited by the evaluation of the nonsmooth max-term. To solve this problem we have introduced a novel adaptive RB-DEIM approach based on a modified version of the space-time a-posteriori error estimator, which suggests that the total error estimator can be decomposed into a RB and a DEIM part. Again numerical results show the capabilities, but also possible limitations of this approach. As long as the DEIM approximations complexity is in the same regime as the RB approximations complexity, the novel approach works well. Especially the speed-up compared to RB approaches can be significantly increased, whilst maintaining the same approximation quality. Compared to classical RB-DEIM, this is only possible due to the adaptive nature of our algorithm.


\medskip
\noindent
\textbf{Acknowledgments.} We are grateful to Denis Korolev (Berlin), Dominik Meidner (Munich) and Karsten Urban (Ulm) for fruitful discussions and very helpful remarks.


\medskip
\noindent
\textbf{Funding.} This research was supported by the German Research Foundation (DFG) under grant number VO 1658/5-2 within the DFG Priority Program ``Non-smooth and Complementarity-based Distributed Parameter Systems: Simulation and Hierarchical Optimization'' (SPP 1962).

\bibliography{BernreutherVolkwein}


\begin{thebibliography}{33}
\ifx \bisbn   \undefined \def \bisbn  #1{ISBN #1}\fi
\ifx \binits  \undefined \def \binits#1{#1}\fi
\ifx \bauthor  \undefined \def \bauthor#1{#1}\fi
\ifx \batitle  \undefined \def \batitle#1{#1}\fi
\ifx \bjtitle  \undefined \def \bjtitle#1{#1}\fi
\ifx \bvolume  \undefined \def \bvolume#1{\textbf{#1}}\fi
\ifx \byear  \undefined \def \byear#1{#1}\fi
\ifx \bissue  \undefined \def \bissue#1{#1}\fi
\ifx \bfpage  \undefined \def \bfpage#1{#1}\fi
\ifx \blpage  \undefined \def \blpage #1{#1}\fi
\ifx \burl  \undefined \def \burl#1{\textsf{#1}}\fi
\ifx \doiurl  \undefined \def \doiurl#1{\url{https://doi.org/#1}}\fi
\ifx \betal  \undefined \def \betal{\textit{et al.}}\fi
\ifx \binstitute  \undefined \def \binstitute#1{#1}\fi
\ifx \binstitutionaled  \undefined \def \binstitutionaled#1{#1}\fi
\ifx \bctitle  \undefined \def \bctitle#1{#1}\fi
\ifx \beditor  \undefined \def \beditor#1{#1}\fi
\ifx \bpublisher  \undefined \def \bpublisher#1{#1}\fi
\ifx \bbtitle  \undefined \def \bbtitle#1{#1}\fi
\ifx \bedition  \undefined \def \bedition#1{#1}\fi
\ifx \bseriesno  \undefined \def \bseriesno#1{#1}\fi
\ifx \blocation  \undefined \def \blocation#1{#1}\fi
\ifx \bsertitle  \undefined \def \bsertitle#1{#1}\fi
\ifx \bsnm \undefined \def \bsnm#1{#1}\fi
\ifx \bsuffix \undefined \def \bsuffix#1{#1}\fi
\ifx \bparticle \undefined \def \bparticle#1{#1}\fi
\ifx \barticle \undefined \def \barticle#1{#1}\fi
\bibcommenthead
\ifx \bconfdate \undefined \def \bconfdate #1{#1}\fi
\ifx \botherref \undefined \def \botherref #1{#1}\fi
\ifx \url \undefined \def \url#1{\textsf{#1}}\fi
\ifx \bchapter \undefined \def \bchapter#1{#1}\fi
\ifx \bbook \undefined \def \bbook#1{#1}\fi
\ifx \bcomment \undefined \def \bcomment#1{#1}\fi
\ifx \oauthor \undefined \def \oauthor#1{#1}\fi
\ifx \citeauthoryear \undefined \def \citeauthoryear#1{#1}\fi
\ifx \endbibitem  \undefined \def \endbibitem {}\fi
\ifx \bconflocation  \undefined \def \bconflocation#1{#1}\fi
\ifx \arxivurl  \undefined \def \arxivurl#1{\textsf{#1}}\fi
\csname PreBibitemsHook\endcsname

\bibitem{Mei11}
\begin{barticle}
\bauthor{\bsnm{Meidner}, \binits{D.}},
\bauthor{\bsnm{Vexler}, \binits{B.}}:
\batitle{A priori error analysis of the {P}etrov–{G}alerkin
  {C}rank–{N}icolson scheme for parabolic optimal control problems}.
\bjtitle{SIAM J. Control and Optimization}
\bvolume{49},
\bfpage{2183}--\blpage{2211}
(\byear{2011}).
\doiurl{10.1137/100809611}
\end{barticle}
\endbibitem

\bibitem{GK11}
\begin{barticle}
\bauthor{\bsnm{Gunzburger}, \binits{M.}},
\bauthor{\bsnm{Kunoth}, \binits{A.}}:
\batitle{Space-time adaptive wavelet methods for optimal control problems
  constrained by parabolic evolution equations}.
\bjtitle{Journal on Control and Optimization}
\bvolume{55},
\bfpage{1150}--\blpage{1170}
(\byear{2011}).
\doiurl{10.1137/100806382}
\end{barticle}
\endbibitem

\bibitem{NV12}
\begin{barticle}
\bauthor{\bsnm{Neitzel}, \binits{I.}},
\bauthor{\bsnm{Vexler}, \binits{B.}}:
\batitle{A priori error estimates for space–time finite element
  discretization of semilinear parabolic optimal control problems}.
\bjtitle{Numerische Mathematik}
\bvolume{120},
\bfpage{345}--\blpage{386}
(\byear{2012}).
\doiurl{10.1007/s00211-011-0409-9}
\end{barticle}
\endbibitem

\bibitem{LS15}
\begin{bbook}
\bauthor{\bsnm{Langer}, \binits{U.}},
\bauthor{\bsnm{Steinbach}, \binits{O.}}:
\bbtitle{Space Time Methods: Applications to Partial Differential Equations}.
\bsertitle{Radon Series on Computational and Applied Mathematics},
vol. \bseriesno{25}.
\bpublisher{De Gruyter},
\blocation{Berlin}
(\byear{2019})
\end{bbook}
\endbibitem

\bibitem{Ste15}
\begin{barticle}
\bauthor{\bsnm{Steinbach}, \binits{O.}}:
\batitle{Space-time finite element methods for parabolic problems}.
\bjtitle{Computational Methods in Applied Mathematics}
\bvolume{15},
\bfpage{551}--\blpage{566}
(\byear{2015}).
\doiurl{10.1515/cmam-2015-0026}
\end{barticle}
\endbibitem

\bibitem{SY18}
\begin{barticle}
\bauthor{\bsnm{Steinbach}, \binits{O.}},
\bauthor{\bsnm{Yang}, \binits{H.}}:
\batitle{Comparison of algebraic multigrid methods for an adaptive space–time
  finite element discretization of the heat equation in 3d and 4d}.
\bjtitle{Numerical Linear Algebra with Applications}
\bvolume{25},
\bfpage{2143}
(\byear{2018}).
\doiurl{10.1002/nla.2143}
\end{barticle}
\endbibitem

\bibitem{HT18}
\begin{barticle}
\bauthor{\bsnm{Harbrecht}, \binits{H.}},
\bauthor{\bsnm{Tausch}, \binits{J.}}:
\batitle{A fast sparse grid based space–time boundary element method for the
  nonstationary heat equation}.
\bjtitle{Numerische Mathematik}
\bvolume{140},
\bfpage{239}--\blpage{264}
(\byear{2018}).
\doiurl{10.1007/s00211-018-0963-5}
\end{barticle}
\endbibitem

\bibitem{Hin20}
\begin{barticle}
\bauthor{\bsnm{Hinze}, \binits{M.}},
\bauthor{\bsnm{Korolev}, \binits{D.}}:
\batitle{A space-time certified reduced basis method for quasilinear parabolic
  partial differential equations}.
\bjtitle{Advances in Computational Mathematics}
\bvolume{47},
\bfpage{36}
(\byear{2021}).
\doiurl{10.1007/s10444-021-09860-z}
\end{barticle}
\endbibitem

\bibitem{SU12}
\begin{barticle}
\bauthor{\bsnm{Steih}, \binits{K.}},
\bauthor{\bsnm{Urban}, \binits{K.}}:
\batitle{Space-time reduced basis methods for time-periodic partial
  differential equations}.
\bjtitle{IFAC Proceedings Volumes}
\bvolume{45},
\bfpage{710}--\blpage{715}
(\byear{2012}).
\doiurl{10.3182/20120215-3-AT-3016.00126}
\end{barticle}
\endbibitem

\bibitem{YPU14}
\begin{barticle}
\bauthor{\bsnm{Yano}, \binits{M.}},
\bauthor{\bsnm{Patera}, \binits{A.T.}},
\bauthor{\bsnm{Urban}, \binits{K.}}:
\batitle{A space-time hp-interpolation-based certified reduced basis method for
  {B}urgers’ equation}.
\bjtitle{Mathematical Models and Methods in Applied Sciences}
\bvolume{24},
\bfpage{1903}--\blpage{1935}
(\byear{2014}).
\doiurl{10.1142/S0218202514500110}
\end{barticle}
\endbibitem

\bibitem{Urb14}
\begin{barticle}
\bauthor{\bsnm{Urban}, \binits{K.}},
\bauthor{\bsnm{Patera}, \binits{A.T.}}:
\batitle{An improved error bound for reduced basis approximation of linear
  parabolic problems}.
\bjtitle{Mathematics of Computation}
\bvolume{83},
\bfpage{1599}--\blpage{1615}
(\byear{2014}).
\doiurl{10.1090/S0025-5718-2013-02782-2}
\end{barticle}
\endbibitem

\bibitem{HPSU22}
\begin{barticle}
\bauthor{\bsnm{Henning}, \binits{J.}},
\bauthor{\bsnm{Palitta}, \binits{D.}},
\bauthor{\bsnm{Simoncini}, \binits{V.}},
\bauthor{\bsnm{Urban}, \binits{K.}}:
\batitle{An ultraweak space-time variational formulation for the wave equation:
  Analysis and efficient numerical solution}.
\bjtitle{ESAIM Mathematical Modelling and Numerical Analysis}
\bvolume{56},
\bfpage{1173}--\blpage{1198}
(\byear{2022}).
\doiurl{10.1051/m2an/2022035}
\end{barticle}
\endbibitem

\bibitem{BRU22}
\begin{botherref}
\oauthor{\bsnm{Beranek}, \binits{N.}},
\oauthor{\bsnm{Reinhold}, \binits{A.}},
\oauthor{\bsnm{Urban}, \binits{K.}}:
A space-time variational method for optimal control problems: Well-posedness,
  stability and numerical solution.
arXiv.
Submitted
(2022).
\doiurl{10.48550/arXiv.2010.00345}
\end{botherref}
\endbibitem

\bibitem{BMNP04}
\begin{barticle}
\bauthor{\bsnm{Barrault}, \binits{M.}},
\bauthor{\bsnm{Maday}, \binits{Y.}},
\bauthor{\bsnm{Nguyen}, \binits{N.C.}},
\bauthor{\bsnm{Patera}, \binits{A.T.}}:
\batitle{An empirical interpolation method: application to efficient
  reduced-basis discretization of partial differential equations}.
\bjtitle{Comptes Rendus Mathematique}
\bvolume{339},
\bfpage{667}--\blpage{672}
(\byear{2004}).
\doiurl{10.1016/j.crma.2004.08.006}
\end{barticle}
\endbibitem

\bibitem{CS10}
\begin{barticle}
\bauthor{\bsnm{Chaturantabut}, \binits{S.}},
\bauthor{\bsnm{Sorensen}, \binits{D.C.}}:
\batitle{Nonlinear model reduction via discrete empirical interpolation}.
\bjtitle{SIAM Journal on Scientific Computing}
\bvolume{32},
\bfpage{2737}--\blpage{2764}
(\byear{2010}).
\doiurl{10.1137/090766498}
\end{barticle}
\endbibitem

\bibitem{CS12}
\begin{barticle}
\bauthor{\bsnm{Chaturantabut}, \binits{S.}},
\bauthor{\bsnm{Sorensen}, \binits{D.C.}}:
\batitle{A state space estimate for {POD}-{DEIM} nonlinear model reduction}.
\bjtitle{SIAM Journal on Numerical Analysis}
\bvolume{50},
\bfpage{46}--\blpage{63}
(\byear{2012}).
\doiurl{10.1137/110822724}
\end{barticle}
\endbibitem

\bibitem{Bet19}
\begin{barticle}
\bauthor{\bsnm{Betz}, \binits{L.M.}}:
\batitle{Second-order sufficient optimality conditions for optimal control of
  non-smooth, semilinear parabolic equations}.
\bjtitle{Journal on Control and Optimization}
\bvolume{57},
\bfpage{4033}--\blpage{4062}
(\byear{2019}).
\doiurl{10.1137/19M1239106}
\end{barticle}
\endbibitem

\bibitem{MS17}
\begin{barticle}
\bauthor{\bsnm{Meyer}, \binits{C.}},
\bauthor{\bsnm{Susu}, \binits{L.M.}}:
\batitle{Optimal control of nonsmooth, semilinear parabolic equations}.
\bjtitle{Journal on Control and Optimization}
\bvolume{55},
\bfpage{2206}--\blpage{2234}
(\byear{2017}).
\doiurl{10.1137/15M1040426}
\end{barticle}
\endbibitem

\bibitem{DHO10}
\begin{bchapter}
\bauthor{\bsnm{Drohmann}, \binits{M.}},
\bauthor{\bsnm{Haasdonk}, \binits{B.}},
\bauthor{\bsnm{Ohlberger}, \binits{M.}}:
\bctitle{Adaptive reduced basis methods for nonlinear convection–diffusion
  equations}.
In: \bbtitle{Finite Volumes for Complex Applications VI Problems \&
  Perspectives}
(\byear{2010}).
\doiurl{10.1007/978-3-642-20671-9_39}
\end{bchapter}
\endbibitem

\bibitem{Ber20}
\begin{bchapter}
\bauthor{\bsnm{Bernreuther}, \binits{M.}},
\bauthor{\bsnm{M\"uller}, \binits{G.}},
\bauthor{\bsnm{Volkwein}, \binits{S.}}:
\bctitle{Reduced basis model order reduction in optimal control of a nonsmooth
  semilinear elliptic {PDE}}.
In: \beditor{\bsnm{Herzog}, \binits{R.}},
\beditor{\bsnm{Heinkenschloss}, \binits{M.}},
\beditor{\bsnm{Kalise}, \binits{D.}},
\beditor{\bsnm{Stadler}, \binits{G.}},
\beditor{\bsnm{Trélat}, \binits{E.}} (eds.)
\bbtitle{Optimization and Control for Partial Differential Equations},
pp. \bfpage{1}--\blpage{32}.
\bpublisher{De Gruyter},
\blocation{Berlin, Boston}
(\byear{2022}).
\doiurl{10.1515/9783110695984-001}
\end{bchapter}
\endbibitem

\bibitem{Eva10}
\begin{bbook}
\bauthor{\bsnm{Evans}, \binits{L.C.}}:
\bbtitle{Partial Differential Equations}.
\bsertitle{Graduate Studies in Mathematics}.
\bpublisher{American Mathematical Society},
\blocation{Providence, Rhode Island}
(\byear{2010})
\end{bbook}
\endbibitem

\bibitem{Zei89a}
\begin{bbook}
\bauthor{\bsnm{Zeidler}, \binits{E.}}:
\bbtitle{Nonlinear Functional Analysis and Its Applications. Linear Monotone
  Operators}
vol. \bseriesno{II/A}.
\bpublisher{Springer},
\blocation{New York}
(\byear{1989})
\end{bbook}
\endbibitem

\bibitem{Zei89b}
\begin{bbook}
\bauthor{\bsnm{Zeidler}, \binits{E.}}:
\bbtitle{Nonlinear Functional Analysis and Its Applications. Nonlinear Monotone
  Operators}
vol. \bseriesno{II/B}.
\bpublisher{Springer},
\blocation{New York}
(\byear{1989})
\end{bbook}
\endbibitem

\bibitem{Tho97}
\begin{bbook}
\bauthor{\bsnm{Thom{\'e}e}, \binits{V.}}:
\bbtitle{Galerkin Finite Element Methods for Parabolic Problems}.
\bpublisher{Springer},
\blocation{Berlin}
(\byear{1997})
\end{bbook}
\endbibitem

\bibitem{Hin10}
\begin{botherref}
\oauthor{\bsnm{Hinterm{\"u}ller}, \binits{M.}}:
Semismooth {N}ewton methods and applications.
Oberwolfach-Seminar on \emph{Mathematics of PDE-Constrained Optimization} at
  Mathematisches Forschungsinstitut in Oberwolfach
(2010).
\url{https://www.math.uni-hamburg.de/home/hinze/Psfiles/
  Hintermueller_OWNotes.pdf}
\end{botherref}
\endbibitem

\bibitem{Ber19}
\begin{botherref}
\oauthor{\bsnm{Bernreuther}, \binits{M.}}:
{RB}-based {PDE}-constrained non-smooth optimization.
Master's thesis,
Universit\"at Konstanz
(2019).
\url{http://nbn-resolving.de/urn:nbn:de:bsz:352-2-t4k1djyj77yn3}
\end{botherref}
\endbibitem

\bibitem{Vol17}
\begin{bchapter}
\bauthor{\bsnm{Gubisch}, \binits{M.}},
\bauthor{\bsnm{Volkwein}, \binits{S.}}:
\bctitle{Chapter 1: {POD} for linear-quadratic optimal control}.
In: \bbtitle{Model Reduction and Approximation - Theory and Algorithms}.
\bsertitle{Computational Science \& Engineering},
pp. \bfpage{3}--\blpage{63}.
\bpublisher{SIAM},
\blocation{Philadelphia}
(\byear{2017}).
\doiurl{10.1137/1.9781611974829.ch1}
\end{bchapter}
\endbibitem

\bibitem{Hesthaven}
\begin{bbook}
\bauthor{\bsnm{Hesthaven}, \binits{J.S.}},
\bauthor{\bsnm{Rozza}, \binits{G.}},
\bauthor{\bsnm{Stamm}, \binits{B.}}:
\bbtitle{Certified Reduced Basis Methods for Parametrized Partial Differential
  Equations}.
\bsertitle{SpringerBriefs in Mathematics}.
\bpublisher{Springer},
\blocation{Cham}
(\byear{2016})
\end{bbook}
\endbibitem

\bibitem{Z12}
\begin{barticle}
\bauthor{\bsnm{Zeng}, \binits{J.-p.}},
\bauthor{\bsnm{Yu}, \binits{H.-x.}}:
\batitle{Error estimates of the lumped mass finite element method for
  semilinear elliptic problems}.
\bjtitle{Journal of Computational and Applied Mathematics}
\bvolume{236},
\bfpage{1993}--\blpage{2004}
(\byear{2012}).
\doiurl{10.1016/j.cam.2011.11.009}
\end{barticle}
\endbibitem

\bibitem{Che20}
\begin{barticle}
\bauthor{\bsnm{Chellappa}, \binits{S.}},
\bauthor{\bsnm{Feng}, \binits{L.}},
\bauthor{\bsnm{Benner}, \binits{P.}}:
\batitle{Adaptive basis construction and improved error estimation for
  parametric nonlinear dynamical systems}.
\bjtitle{International Journal for Numerical Methods in Engineering}
\bvolume{121},
\bfpage{5320}--\blpage{5349}
(\byear{2020}).
\doiurl{10.1002/nme.6462}
\end{barticle}
\endbibitem

\bibitem{Pet15}
\begin{barticle}
\bauthor{\bsnm{Peherstorfer}, \binits{B.}},
\bauthor{\bsnm{Willcox}, \binits{K.}}:
\batitle{Online adaptive model reduction for nonlinear systems via low-rank
  updates}.
\bjtitle{SIAM Journal on Scientific Computing}
\bvolume{37},
\bfpage{2123}--\blpage{2150}
(\byear{2015}).
\doiurl{10.1137/140989169}
\end{barticle}
\endbibitem

\bibitem{FENICS}
\begin{botherref}
\oauthor{\bsnm{Alna\hspace{-0.45mm}es}, \binits{M.}},
\oauthor{\bsnm{Blechta}, \binits{J.}},
\oauthor{\bsnm{Hake}, \binits{J.}},
\oauthor{\bsnm{Johansson}, \binits{A.}},
\oauthor{\bsnm{Kehlet}, \binits{B.}},
\oauthor{\bsnm{Logg}, \binits{A.}},
\oauthor{\bsnm{Richardson}, \binits{C.}},
\oauthor{\bsnm{Ring}, \binits{J.}},
\oauthor{\bsnm{Rognes}, \binits{M.E.}},
\oauthor{\bsnm{Wells}, \binits{G.N.}}:
The {FEniCS} project version 1.5.
Archive of Numerical Software
\textbf{3}
(2015).
\doiurl{10.11588/ans.2015.100.20553}
\end{botherref}
\endbibitem

\bibitem{Scipy}
\begin{barticle}
\bauthor{\bsnm{{Virtanen}}, \binits{P.}},
\bauthor{\bsnm{{Gommers}}, \binits{R.}},
\bauthor{\bsnm{{Oliphant}}, \binits{T.E.}},
\bauthor{\bsnm{{Haberland}}, \binits{M.}},
\bauthor{\bsnm{{Reddy}}, \binits{T.}},
\bauthor{\bsnm{{Cournapeau}}, \binits{D.}},
\bauthor{\bsnm{{Burovski}}, \binits{E.}},
\bauthor{\bsnm{{Peterson}}, \binits{P.}},
\bauthor{\bsnm{{Weckesser}}, \binits{W.}},
\bauthor{\bsnm{{Bright}}, \binits{J.}},
\bauthor{\bsnm{{van der Walt}}, \binits{S.J.}},
\bauthor{\bsnm{{Brett}}, \binits{M.}},
\bauthor{\bsnm{{Wilson}}, \binits{J.}},
\bauthor{\bsnm{{Jarrod Millman}}, \binits{K.}},
\bauthor{\bsnm{{Mayorov}}, \binits{N.}},
\bauthor{\bsnm{{Nelson}}, \binits{A.R.J.}},
\bauthor{\bsnm{{Jones}}, \binits{E.}},
\bauthor{\bsnm{{Kern}}, \binits{R.}},
\bauthor{\bsnm{{Larson}}, \binits{E.}},
\bauthor{\bsnm{{Carey}}, \binits{C.J.}},
\bauthor{\bsnm{{Polat}}, \binits{{\. I}.}},
\bauthor{\bsnm{{Feng}}, \binits{Y.}},
\bauthor{\bsnm{{Moore}}, \binits{E.W.}},
\bauthor{\bsnm{{VanderPlas}}, \binits{J.}},
\bauthor{\bsnm{{Laxalde}}, \binits{D.}},
\bauthor{\bsnm{{Perktold}}, \binits{J.}},
\bauthor{\bsnm{{Cimrman}}, \binits{R.}},
\bauthor{\bsnm{{Henriksen}}, \binits{I.}},
\bauthor{\bsnm{{Quintero}}, \binits{E.A.}},
\bauthor{\bsnm{{Harris}}, \binits{C.R.}},
\bauthor{\bsnm{{Archibald}}, \binits{A.M.}},
\bauthor{\bsnm{{Ribeiro}}, \binits{A.H.}},
\bauthor{\bsnm{{Pedregosa}}, \binits{F.}},
\bauthor{\bsnm{{van Mulbregt}}, \binits{P.}},
\bauthor{\bsnm{{Contributors}}, \binits{S...}}:
\batitle{{{SciPy} 1.0: Fundamental Algorithms for Scientific Computing in
  Python}}.
\bjtitle{Nature Methods}
\bvolume{17},
\bfpage{261}--\blpage{272}
(\byear{2020}).
\doiurl{10.1038/s41592-019-0686-2}
\end{barticle}
\endbibitem

\end{thebibliography}

\end{document}